\documentstyle{amltd}
\begin{document}
% \input BoxedEps.Tex %for mac
%\SetTexturesEPSFSpecial %for mac
\input boxedeps.tex % unix
\SetepsfEPSFSpecial % unix
\HideDisplacementBoxes
\def\figin#1#2{\medbreak
$$
 {\BoxedEPSF{#1.eps scaled
#2}%
}%
$$
\medbreak\noindent}
\annalsline{153}{2001}
\received{July 5, 1999}
\startingpage{355}
\def\bye{ \end{document}}
 \font\tenrm=cmr10
\def\ritem#1{\item[{\rm #1}]}
\def\nn{\nonumber}
\def\nne{\nonumber\end{eqnarray}}
\def\bel#1{\begin{eqnarray}\label{#1}}
%--------------- Author macros ---------------
\def\et{\endproclaim}
\def\ec{\endproclaim}
\def\ep{\endproclaim}
\def\el{\endproclaim}
%for Bbb in amstex
\catcode`\@=11
\font\twelvemsb=msbm10 scaled 1100
\font\tenmsb=msbm10
%\font\ninemsb=msbm7 scaled 1100%msbm9
\font\ninemsb=msbm10 scaled 800
\newfam\msbfam
\textfont\msbfam=\twelvemsb  \scriptfont\msbfam=\ninemsb
  \scriptscriptfont\msbfam=\ninemsb
\def\msb@{\hexnumber@\msbfam}
\def\Bbb{\relax\ifmmode\let\next\Bbb@\else
 \def\next{\errmessage{Use \string\Bbb\space only in math
mode}}\fi\next}
\def\Bbb@#1{{\Bbb@@{#1}}}
\def\Bbb@@#1{\fam\msbfam#1}
\catcode`\@=12

 \catcode`\@=11
\font\twelveeuf=eufm10 scaled 1100
\font\teneuf=eufm10
\font\nineeuf=eufm7 scaled 1100%eufm9
\newfam\euffam
\textfont\euffam=\twelveeuf  \scriptfont\euffam=\teneuf
  \scriptscriptfont\euffam=\nineeuf
\def\euf@{\hexnumber@\euffam}
\def\frak{\relax\ifmmode\let\next\frak@\else
 \def\next{\errmessage{Use \string\frak\space only in math
mode}}\fi\next}
\def\frak@#1{{\frak@@{#1}}}
\def\frak@@#1{\fam\euffam#1}
\catcode`\@=12

 \def\be#1{\begin{equation}\label{#1}}
\newcommand{\ee}{\end{equation}}

\newcommand{\bi}{\begin{itemize}}
\newcommand{\ei}{\end{itemize}}

\newcommand{\ben}{\begin{enumerate}}
\newcommand{\een}{\end{enumerate}}

%%%%%%%%%%%% Definitions %%%%%%%%%%%%%%%%%%%%%%%%%%%%%%%

\def \ba {\begin{array}}
\def \ea {\end{array}}
\def \beqa {\begin{eqnarray}}
\def \eeqa {\end{eqnarray}}

 \def \Z {{\Bbb Z}}
\def \R {{\Bbb R}}
\def \C {{\Bbb C}}
\def \N {{\Bbb N}}
\def \X {{\Bbb X}}
\def \E {{\Bbb E}\,}
\def \P {{\Bbb P}}
\def \ra {\rightarrow}
\def \da {\downarrow}
\def \ua {\uparrow}

\def \cL {{\cal L}}
\def \cG {{\cal G}}
\def \cP {{\cal P}}
\def \cU {{\cal U}}
\def \cR {{\cal R}}
\def \cM {{\cal M}}
\def \cB {{\cal B}}
\def \cR {{\cal R}}

\def \RSNI {\hbox{\small {\rm RSNI}}}

%%%%%%%%%%%%%%%%%%%%%%%%%%%%%%%%%%%%%%%%%%%%%%%%%%%%%%%%%%%%%%%%%%%%%%%%%%%%%%%%
%-------------- Author entries --------------------

\title{Moderate deviations for the volume\\ of the Wiener sausage} %Article title
\shorttitle{Moderate deviations of the Wiener sausage}  

%% \acknowledgements{}
%\author{}
 \twoauthors{M.\ van den Berg, E.\ Bolthausen,}{F.\ den Hollander}
 \institutions{
School of Mathematics, University of Bristol, Bristol, United Kingdom\\
{\eightpoint {\it E-mail address\/}:  m.vandenberg@bristol.ac.uk}\\
\vglue6pt
 Institut f\"ur  Mathematik, Universit\"at Z\"urich, Z\"urich,
Switzerland\\
{\eightpoint {\it E-mail address\/}:  eb@math.unizh.ch}\\
\vglue6pt
EURANDOM, P.O. Box 513, 5600 MB, Eindhoven, The
Netherlands\\
{\eightpoint {\it E-mail address\/}:   denhollander@eurandom.tue.nl
}} 
 
 \vfill \centerline{\bf Abstract}
\bigbreak
For $a>0$, let $W^a(t)$ be the $a$--neighbourhood of standard Brownian\break
motion in $\R^d$ starting at 0 and observed until time $t$. It is well-known that $E|W^a(t)| \sim \kappa_a t$ ($t \ra
\infty$) for $d \geq 3$, with $\kappa_a$ the Newtonian capacity of the ball with radius $a$. We
prove that
$$
\lim_{t \ra \infty} \frac{1}{t^{(d-2)/d}} \log P(|W^a(t)| \leq bt) =
- I^{\kappa_a}(b) \in (-\infty,0) \qquad \mbox{ for all } 0<b<\kappa_a
$$
and derive a variational representation for the rate function
$I^{\kappa_a}$. We
show that the optimal strategy to realise the above moderate deviation
is for $W^a(t)$ to `look like a Swiss cheese': $W^a(t)$ has random holes
whose sizes are of order 1 and whose density varies on scale $t^{1/d}$.
The optimal strategy is such that $t^{-1/d}W^a(t)$ is delocalised in the
limit as $t \ra \infty$. This is markedly different from the optimal
strategy for large deviations $\{|W^a(t)| \leq f(t)\}$ with $f(t)=o(t)$,
where $W^a(t)$ is known to fill completely  a ball of volume $f(t)$ and
nothing outside, so that $W^a(t)$ has no holes and $f(t)^{-1/d}W^a(t)$
is localised in the limit as $t \ra \infty$.

We give a detailed analysis of the rate function $I^{\kappa_a}$, in
particular, its
behaviour near the boundary points of $(0,\kappa_a)$ as well as certain
monotonicity properties. It turns out that $I^{\kappa_a}$ has an infinite
slope at
$\kappa_a$ and, remarkably, for $d \geq 5$ is nonanalytic at some critical
point in $(0,\kappa_a)$, above which it follows a pure power law. This
crossover is associated with a collapse transition in the optimal strategy.

We also derive the analogous moderate deviation result for $d=2$. In this case
$E|W^a(t)|\sim 2\pi t/\log t$ ($t\ra\infty$), and we prove that
$$
\lim_{t \ra \infty} \frac{1}{\log t} \log P(|W^a(t)|\leq bt/\log t)
=-I^{2\pi}(b)\in (-\infty,0) \qquad \mbox{ for all } 0<b<2\pi.
$$
The rate function $I^{2\pi}$ has a finite slope at $2\pi$.
\eject
\demo{Acknowledgment} Part of this research was supported by the
Volkswagen-Stiftung
through the RiP-program at the Mathematisches Forschungsinstitut,
Oberwolfach, Germany.
MvdB was supported by the London Mathematical Society, EB was supported by the
Swiss National Science Foundation through grant 20--49501.96, MvdB and FdH
were supported by the British Council and The Netherlands Organisation for
Scientific
Research through grant JRP431.  

\section{Introduction and main results:\\ Theorems \ref{t1}--\ref{t5} and
Corollaries \ref{c1}, \ref{c2}}

 1.1 {\it The Wiener sausage}.
Let $\beta(t)$, $t \geq 0$, be the standard Brownian motion in $\R^d$ --
the Markov process with generator $\Delta/2$ -- starting at 0. Let $P,E$ denote
its probability law and expectation on path space. The Wiener sausage with
radius
$a>0$ is the process defined by
\be{WS}
W^a(t) = \bigcup_{0 \leq s \leq t} B_a(\beta(s)), ~~~t \geq 0,
\ee
where $B_a(x)$ is the open ball with radius $a$ around $x \in \R^d$.
The Wiener sausage is an important mathematical object, because it is one of
the simplest examples of a non-Markovian functional of Brownian motion. It
plays
a key role in the study of various stochastic phenomena, such as heat
conduction and
trapping in random media, as well as in the analysis of spectral properties
of random
Schr\"odinger operators.

A lot is known about the behaviour of the volume of $W^a(t)$ as $t \ra
\infty$. For instance,
\be{mean}
E|W^a(t)| \sim
\left\{\ba{ll}
\sqrt{8t/\pi}    &(d=1)\\
2\pi t/\log t    &(d=2)\\
\kappa_a t       &(d \geq 3),
\ea\right.
\ee
with $\kappa_a = a^{d-2}2\pi^{d/2}/\Gamma(\frac{d-2}{2})$ the Newtonian
capacity of
$B_a(0)$ associated with the Green function of $(-\Delta/2)^{-1}$, and
\be{variance}
{\rm Var}|W^a(t)| \asymp
\left\{\ba{ll}
t             &(d=1)\\
t^2/\log^4 t  &(d=2)\\
t \log t      &(d=3)\\
t             &(d \geq 4)
\ea\right.
\ee
(Spitzer \cite{Sp}, Le Gall \cite{LG1}). Moreover, $|W^a(t)|$ satisfies the
strong law and the central limit theorem for $d \geq 2$; the limit law is
Gaussian
for $d \geq 3$ and non-Gaussian for $d=2$ (Le Gall \cite{LG2}).
Note that for $d \geq 2$ the Wiener sausage is a {\it sparse} object: since
the Brownian
motion typically travels a distance $\sqrt{t}$ in each direction, (\ref{mean})
shows that most of the space in the convex hull of $W^a(t)$ is not covered.

\demo{{\rm 1.2.} Large deviations}
The large deviation properties of  $|W^a(t)|$ in\break the {\it downward} direction
have been studied by Donsker and Varadhan \cite{DV},\break Bolthausen \cite{Bo} and
Sznitman \cite{Sz}. For $d \geq 2$ the outcome, proved in successive stages of
refinement, reads as follows:
\be{lds}
\lim\limits_{t \ra \infty} \frac{f(t)^{2/d}}{t} \log P(|W^a(t)| \leq f(t))
= - \frac{1}{2} \lambda_d
\ee
for any $f \colon \R_+ \mapsto \R_+$ satisfying $\lim\limits_{t \ra \infty}
f(t) = \infty$ and
\be{fcond}
f(t) =
\left\{\ba{ll}
o(t/\log t) &(d=2)\\
o(t)        &(d \geq 3),
\ea\right.
\ee
where $\lambda_d>0$ is the smallest Dirichlet eigenvalue of $-\Delta$ on
the ball with unit volume. It turns out that the optimal strategy for the
Brownian motion to realise the large deviation in (\ref{lds}) is to stay
inside a ball with volume $f(t)$ until time $t$, i.e., the Wiener sausage
covers this ball entirely and nothing outside. (The optimality comes from the
Faber-Krahn isoperimetric inequality, and the cost of staying inside the
ball is
$\exp[- \frac{1}{2} \lambda_d t/f(t)^{2/d}]$ to leading order.) Thus, the
optimal
strategy is simple and $f(t)^{-1/d}W^a(t)$ is {\it localised}. Note that,
apparently,
a large deviation {\it below} the scale of the mean `squeezes all the empty
space out of the Wiener sausage'. Also note that the limit in (\ref{lds})
does not depend on $a$.

The law of the Brownian motion {\it conditioned} on the large deviation
event $\{|W^a(t)| \leq f(t)\}$ has been studied by Sznitman \cite{Szconf},
Bolthausen \cite{Boconf} and Povel \cite{Poconf}. This law is indeed like
the optimal strategy described above, with an explicitly known probability
distribution for the centre of the ball the Brownian motion stays confined in.
\enddemo

\demo{{\rm 1.3.} Moderate deviations}
The aim of the present paper is to extend (\ref{lds}),  (\ref{fcond}) by
investigating deviations {\it on} the scale of the mean. We call such
deviations moderate.\footnote{The term `moderate' is often used for deviations away from the
mean that are smaller than the scale of the mean, but in view of the contrast
with (\ref{lds}),  (\ref{fcond}) we prefer this terminology.}
Our first main result reads:
\specialnumber{1}
\proclaim{Theorem} \label{t1}
Let $d \geq 3$ and $a>0${\rm .} For every $b>0${\rm ,}
\be{mds}
\lim_{t \ra \infty} \frac{1}{t^{(d-2)/d}} \log P(|W^a(t)| \leq bt) = -
I^{\kappa_a}(b),
\ee
where
\be{mdsrf}
I^{\kappa_a}(b) = \inf_{\phi \in \Phi^{\kappa_a}(b)} \Big[ \frac{1}{2}
\int_{\R^d} |\nabla
\phi|^2(x) dx
\Big]
\ee
with
\be{mdsset}
\Phi^{\kappa_a}(b) = \Big\{\phi \in H^1(\R^d) \colon \int_{\R^d}\phi^2(x)
dx = 1,~
\int_{\R^d}\Big(1-e^{-\kappa_a\phi^2(x)}\Big)dx \leq b \Big\}.
\ee
\endproclaim 

The idea behind Theorem \ref{t1} is that the optimal strategy for the
Brownian motion
to realise the event $\{|W^a(t)| \leq bt\}$ is to behave like a Brownian
motion in
a drift field $xt^{1/d} \mapsto (\nabla \phi/\phi)(x)$ for some smooth
$\phi \colon \R^d
\mapsto [0,\infty)$. The cost of adopting this drift during a time $t$ is
the exponential of
$t^{(d-2)/d}$ times the integral in (\ref{mdsrf}) to leading order. The
effect of the drift is to push
the Brownian motion towards the origin. Conditioned on adopting the drift,
the Brownian
motion spends time $\phi^2(x)$ per unit volume in the neighbourhood of
$xt^{1/d}$, and it
turns out that the Wiener sausage covers a fraction $1-\exp[-\kappa_a\phi^2(x)]$
of the space in that neighbourhood. The best choice of the drift field is
therefore given by a
minimiser of the variational problem in (\ref{mdsrf}), or by a minimising
sequence.

We thus see that the optimal strategy for the Wiener sausage is to cover
only part of
the space and to leave {\it random holes}\footnote{The motto of this paper: `How a Wiener sausage turns into a Swiss
cheese'.}
whose sizes are of order 1 and whose density varies on scale $t^{1/d}$.
This strategy is
more complicated than for (\ref{lds}) and $t^{-1/d}W^a(t)$ is {\it
delocalised}.
(In Section 5.1 it is shown that all minimisers or minimising sequences of
(\ref{mdsrf})
are strictly positive.) Note that, apparently, a moderate deviation {\it
on} the scale
of the mean `does not squeeze all the empty space out of the Wiener
sausage'. Also note
that the limit in (\ref{mds}) does depend on $a$.\footnote{To prove that the law of the Brownian motion {\it
conditioned} on the moderate deviation event $\{|W^a(t)|\leq bt\}$ actually follows the
optimal
`Swiss cheese strategy' requires substantial extra work. We shall not
address this issue here.
Even though we shall sometimes interpret our results in terms of this
strategy, we
have no pathwise statements to offer.}

It is clear from (\ref{lds}),  (\ref{fcond}) that the case $d=2$ is critical. Our
next main result is the following parallel of Theorem \ref{t1}.

\phantom{cd}

\specialnumber{2} \proclaim{Theorem} \label{t2}
Let $d=2$ and $a>0${\rm .} For every $b>0,$ \vglue-9pt
\be{mds2}
\lim_{t \ra \infty} \frac{1}{\log t} \log P(|W^a(t)| \leq bt/\log t) = -
I^{2\pi}(b),
\ee
where $I^{2\pi}(b)$ is given by the same formulas as in
{\rm (\ref{mdsrf}),  (\ref{mdsset}),} except
that $\kappa_a$ is replaced by $2\pi${\rm .}
\endproclaim
 
\vglue8pt

Theorem \ref{t2} shows that for $d=2$ the moderate deviations have a
polynomially small rather than an exponentially small probability. The
optimal strategy is of the same type, but now the Wiener sausage lives
on scale $\sqrt{t/\log t}$, which is only slightly below the diffusive scale.
Contrary to the case $d \geq 3$, the rate function does not depend on $a$.
This means that the random holes in the Swiss cheese have a typical size and
a typical mutual distance that tend to infinity as $t \ra \infty$,
washing out the dependence on the radius of the Wiener sausage.

There is no result analogous to Theorems \ref{t1}, \ref{t2} for $d=1$, for the
simple reason that the strong law fails (see (\ref{mean}), (\ref{variance})).
The variational problem in (\ref{mdsrf}), (\ref{mdsset}) certainly continues
to make sense for $d=1$, but it does not describe the Wiener sausage: holes
are impossible in $d=1$.

\demo{{\rm 1.4.} The rate function}
We proceed with a closer analysis of (\ref{mdsrf}), (\ref{mdsset}).
First we scale out the $a$--dependence and make some general statements
about the rate function. Recall that $\kappa_a \equiv 2\pi$ for $d=2$.\enddemo

\specialnumber{3} \proclaim{Theorem} \label{t3} Let $d \ge 2$ and $a>0${\rm ,}
\smallbreak
{\rm (i)} For every $b>0,$
\be{scalI}
I^{\kappa_a}(b) = \frac{1}{2\kappa_a^{2/d}} \chi(b/\kappa_a),
\ee
where $\chi \colon (0,\infty) \mapsto [0,\infty)$ is given by
\be{chi}
\chi(u) = \inf \{ \|\nabla \psi \|_2^2 \colon ~\psi \in H^1(\R^d),
~\|\psi\|_2 = 1, ~\textstyle\int_{\R^d} (1-e^{-\psi^2}) \leq u \}.\hskip.25in
\ee
\medbreak
{\rm (ii)} $\chi$ is continuous on $(0,\infty)${\rm ,} strictly decreasing on
$(0,1)${\rm ,} and equal to zero on $[1,\infty)${\rm .} 
\medbreak
{\rm (iii)} $u \mapsto u^{2/d}\chi(u)$ is strictly decreasing on $(0,1)$ and
\be{eq1.12}
\lim_{u\downarrow0}u^{2/d}\chi(u)=\lambda_d
\ee
with $\lambda_d$ as defined below {\rm (\ref{fcond}).}
\endproclaim

Theorem \ref{t3}(iii) shows that the limit $b \da 0$ connects up nicely with
(\ref{lds}), (\ref{fcond}).

Our next two results show that the variational problem in (\ref{chi}) displays
a surprising dimension dependence.
 \figin{fig1}{900}
\begin{center}
 Qualitative picture of $u \mapsto \chi(u)$ for\\
(i) $d=2$; (ii) $d=3,4$; (iii) $d \geq 5$. 
\end{center}

\specialnumber{4} \proclaim{Theorem} \label{t4}
Let $2 \leq d \leq 4${\rm .}
 
{\rm (i)} For every $u\in(0,1)$ the variational problem in {\rm (\ref{chi})} has a
minimiser
that is strictly positive{\rm ,} radially symmetric {\rm (}\/modulo shifts\/{\rm )} and strictly decreasing
in the radial component{\rm .} Any other minimiser is of the same type{\rm .} 
\medbreak

{\rm (ii)} $u \mapsto (1-u)^{-2/d}\chi(u)$ is strictly decreasing on $(0,1)$ and
\be{eq1.13}
\lim_{u\uparrow 1}(1-u)^{-2/d}\chi(u)=2^{2/d}\mu_d,
\ee
where
\be{eq1.14}
\mu_d =
\left\{\ba{ll}
\inf\{\|\nabla\psi\|_2^2 \colon ~\psi \in H^1(\R^d),
~\|\psi\|_2=1,~\|\psi\|_4 = 1 \} &(d=2,3)\\[4pt]
\inf\{\|\nabla\psi\|_2^2 \colon ~\psi \in D^1(\R^4),
~\|\psi\|_4 = 1 \} &(d=4)
\ea\right.
\ee
satisfying $0<\mu_d<\infty${\rm .}
\endproclaim

\specialnumber{5} \proclaim{Theorem} \label{t5}
Let $d\ge5${\rm .} 

{\rm (i)} Define
\bel{eq1.15}
\nu_d &=& \inf\{\|\nabla\psi\|_2^2 \colon ~\psi \in D^1(\R^d),
~\int_{\R^d}(e^{-\psi^2}-1+\psi^2)=1\}  \\
\Sigma &=& \mbox{the set of minimisers}\nn\\
\Sigma^* &=& \mbox{the set of local minimisers}.
\nne
Then $0<\nu_d<\infty$ and $\emptyset \neq \Sigma^* \supseteq \Sigma${\rm .}
Moreover{\rm ,} all
elements of $\Sigma$ are strictly positive{\rm ,} radially symmetric {\rm (}\/modulo
shifts\/{\rm ),} strictly
decreasing in the radial component{\rm ,} and there exists a constant $K_d$ such that \vglue-9pt
\be{eq1.16}
\|\psi\|_2^2 > \frac{d}{d-2} \quad \mbox{ for all } \psi \in \Sigma^*, \qquad
\|\psi\|_2^2 \leq K_d \quad \mbox{ for all } \psi \in \Sigma.\hskip.5in
\ee
\bigbreak
{\rm (ii)} Define $2/d \leq u^*_d \leq u^-_d \leq u^+_d \leq 1-K_d^{-1} < 1$ by \vglue-9pt
\be{defud's}
u^*_d = 1 - [\inf_{\psi\in\Sigma^*}\|\psi\|^2_2]^{-1}, \quad
u^-_d = 1 - [\inf_{\psi\in\Sigma}\|\psi\|^2_2]^{-1}, \quad
u^+_d = 1 - [\sup_{\psi\in\Sigma}\|\psi\|^2_2]^{-1}.
\ee
For every $u\in(0,u^*_d]$ the variational problem in {\rm (\ref{chi})} has a
minimiser
that is strictly positive{\rm ,} radially symmetric {\rm (}\/modulo shifts\/{\rm )} and strictly
decreasing
in the radial component{\rm .} Any other minimiser is of the same type{\rm .} For\break every
$u\in(u^+_d,1)$
the variational problem in {\rm (\ref{chi})} does not have a minimiser{\rm .}\break There
exists a minimising
sequence $(\psi_j)$ such that $\psi_j(\cdot)$ converges weakly to\break
$\psi(\cdot\,(1-u)^{-1/d})$
in $H^1(\R^d)$ as $j \to \infty$ for some $\psi\in\Sigma${\rm .} \medbreak

{\rm (iii)} $u \mapsto (1-u)^{-(d-2)/d}\chi(u)$ is strictly decreasing on
$(0,u^*_d]${\rm ,} nonincreasing
and strictly greater than $\nu_d$ on $(u^*_d,u^-_d)${\rm ,} while
\be{eq1.17}
(1-u)^{-(d-2)/d}\chi(u) \equiv \nu_d \quad \mbox{ for } u \in [u^-_d,1).
\ee
\endproclaim

Note that $\chi$ is nonanalytic at $u^-_d$. Whether or not (\ref{chi}) has
a minimiser for
$u \in (u^*_d,u^+_d]$ and whether or not (\ref{eq1.15}) has just one local
minimiser
both remain open. Possibly $|\Sigma^*|=1$, in which case $u^*_d=u^-_d=u^+_d
=: u_d$, but this seems hard to settle (see \S\S 5.6--5.8).

\demo{{\rm 1.5.} Comments}
To explain the situation in Theorem \ref{t5}, let us insert the scaling
$\psi(\cdot\,(1-u)^{-1/d})$ into (\ref{chi}) to obtain
\begin{eqnarray} \label{chiscal1}
\qquad  && (1-u)^{(d-2)/d} \chi(u)  
 = {\rm inf}\{\|\nabla\psi\|^2_2 \colon ~\psi \in H^1(\R^d),
~\|\psi\|^2_2 = (1-u)^{-1},\quad \\ \noalign{\vskip4pt}
&&\hskip2.2in ~\textstyle\int_{\R^d} (e^{-\psi^2}-1+\psi^2)
\geq 1\}. \nonumber
\end{eqnarray}  
In Section 5.1 it will be shown that the two constraints in
(\ref{chiscal1}) may be replaced by
$\|\psi\|^2_2 \leq (1-u)^{-1}$ and $\int_{\R^d} (e^{-\psi^2}-1+\psi^2) =
1$, after which
we have a variational problem as in (\ref{eq1.15}) but with an upper bound
on $\|\psi\|^2_2$.
Let us now consider the optimistic scenario where $\Sigma^*$ has a unique
element $\psi^*$.
Then $u^*_d = u^-_d = u^+_d =: u_d$ and $\|\psi^*\|^2_2 = (1-u_d)^{-1}$. It
turns out
that for $u \in (0,u_d]$ the variational problem in (\ref{chiscal1}) has a
minimiser because
no $L^2$--mass wants to leak away to infinity (even though this minimiser
has little to do
with $\psi^*$ itself). On the other hand, for $u \in (u_d,1)$ it has no
minimiser, and any
minimising sequence converges weakly to $\psi^*$ by leaking $L^2$--mass. In
the less
optimistic scenario where $|\Sigma^*|>1$, there is no leakage for $u \in
(0,u^*_d]$ and leakage
for $u \in (u^+_d,1)$.

The situation in Theorem \ref{t4} can be explained as follows. It turns out
that for $2 \leq d
\leq 4$ all elements of $\Sigma^*$ have infinite $L^2$--norm, so that
$u^*_d=u^-_d=u^+_d=1$.
Hence for any $u \in (0,1)$ there is no leakage and (\ref{chiscal1}) has a minimiser.

The following points in Theorems \ref{t4} and \ref{t5} are noteworthy:
\begin{itemize}
\item[I.] At $b=\kappa_a$ the rate function has an {\it infinite} slope for
$d \geq 3$
but a {\it finite} slope for $d=2$.
\item[II.] The scaling as $b \ua \kappa_a$ is different for $2 \leq d \leq
4$ and $d \geq 5$.
Apparently a delicate dimension dependence is felt as the deviation becomes
smaller than the mean.
The fact that for $d \geq 5$ there is no minimiser for $u \in (u^+_d,1)$ is
to be interpreted as
saying that the optimal strategy is {\it time-inhomogeneous} in the
following sense. Let us again
pretend that $\Sigma^*$ has a unique element $\psi^*$, and let us put
$\rho(u) = (u-u_d)/(1-u_d)
\in (0,1)$. Then heuristically (recall footnote 3):
\begin{itemize}
\item[(1)] Until time $[1-\rho(u)]t$ the Wiener sausage makes a Swiss
cheese on scale
$t^{1/d}$ parametri\-sed by $\psi^*(\cdot\,(1-u)^{-1/d})$, filling a volume
$\kappa_a[u-\rho(u)]t$.
\item[(2)] After time $[1-\rho(u)]t$ it behaves like a typical Wiener sausage
on scale $\sqrt{t}$, filling an additional volume $\kappa_a \rho(u) t$.
\end{itemize}
Thus, at time $[1-\rho(u)]t$ the optimal strategy undergoes a {\it collapse
transition}
from sub\-diffusive behaviour (scale $t^{1/d}$) to diffusive behaviour
(scale $\sqrt{t}$).
The picture is unclear when $|\Sigma^*|>1$ and $u \in (u^*_d,u^+_d]$.
Still, we expect
some type of collapse transition to occur.
\item[III.] For $d \geq 3$, the scaling of the rate function near
$\kappa_a$ does
not connect up with the central limit theorem. Indeed, if we pick $b=b_t$ with
\end{itemize}

\phantom{time to come home}
\vglue-30pt
\be{md+clt-choice}
b_t t = \left\{\ba{ll}
\kappa_a t  - c \sqrt{t \log t}      &(d=3)\\
\kappa_a t - c t                     &(d \geq 4)
\ea\right.
\ee
\noindent \hangindent=28pt \hangafter=0 for some $c>0$ and recall (\ref{mean}), (\ref{variance}), then we find from
(\ref{scalI}),
(\ref{eq1.13}) and (\ref{eq1.17}) that

\be{md+clt}
I^{\kappa_a}(b_t) t^{(d-2)/d} \ra \infty \quad (t \ra \infty)
\ee
\noindent \hangindent=28pt \hangafter=0 instead of a finite limit. Therefore the moderate deviations are in a sense
{\it anomalous}. For $d=2$, on the other hand, we put

\be{cltextra1}
b_t \frac{t}{\log t} = \frac{2\pi t}{\log t} - c \frac{t}{\log^2 t}
\ee
\noindent \hangindent=28pt \hangafter=0  for some $c>0$ and find that

\be{cltextra2}
\lim_{t \ra \infty} I^{2\pi}(b_t) \log t \enspace \mbox{ exists in } (0,\infty).
\ee
\noindent \hangindent=28pt \hangafter=0 So there is {\it no anomaly} in this case. Incidentally, for $d \geq 3$ the
correction term to
the asymptotic mean is of smaller order than the asymptotic standard deviation,
while for $d=2$ it is of the same order (Spitzer \cite{Sp}, Getoor \cite{Ge}).
For the above argument we may therefore indeed only consider the leading order
terms given by (\ref{mean}), (\ref{variance}).
\medbreak

\noindent
The anomaly for $d \geq 3$ is somewhat surprising. It suggests that the
central limit behaviour is controlled by the {\it local\/} fluctuations of
the Wiener
sausage, while the moderate and large deviations are controlled by the {\it
global\/}
fluctuations.

It remains open whether $I^{2\pi}$ is convex for $d=2$ and whether
$I^{\kappa_a}$ has only one
point of inflection for $d \geq 3$.
\enddemo

\demo{{\rm 1.6.} Negative exponential moments}
We close this introduction with two corollaries. An immediate consequence
of Theorem \ref{t1}
is the following result.\footnote{Sznitman \cite[pp.\ 213--214]{Szbook}  gives a heuristic
derivation of Corollary \ref{c1}
using his method of `enlargement of obstacles'.}
\enddemo

\specialnumber{1} \proclaim{{C}orollary} \label{c1}
Let $d \geq 3$ and $a>0${\rm .} For every $c>0,$
\be{corr1}
\lim\limits_{t \ra \infty} \frac{1}{t^{(d-2)/d}}
\log E\biggl(\exp[-ct^{-2/d}|W^a(t)|]\biggr) = -J^{\kappa_a}(c)
\ee
with
\be{corr2}
J^{\kappa_a}(c) = \inf\limits_{b>0} [bc + I^{\kappa_a}(b)].
\ee
\ec

\vglue9pt

It follows from (\ref{corr2}) that
\be{corr3}
J^{\kappa_a}(c)
\left\{\ba{ll} = \kappa_a c &(0<c \leq c_a^{\ast})\\
               < \kappa_a c &(c > c_a^{\ast})
\ea\right.
\ee
with
\be{corr4}
c_a^{\ast} = \max \{c>0 \colon I^{\kappa_a}(b) \geq c (\kappa_a - b) \mbox{
for all
} 0<b\leq \kappa_a \}.\hskip.75in
\ee
At $c=c_a^{\ast}$, the minimiser of (\ref{corr2}) moves from $b=\kappa_a$
to the
interior of $(0,\kappa_a]$. Heuristically, this corresponds to a {\it collapse
transition} in the optimal strategy for the Brownian motion associated with
(\ref{corr1}), (\ref{corr2}), namely, from diffusive behaviour (scale $\sqrt{t}$)
to subdiffusive behaviour (scale $t^{1/d}$). By Theorems \ref{t3}(i),
\ref{t4}(ii)
and \ref{t5}(iii), the left derivative of $I^{\kappa_a}$ at $b=\kappa_a$ is
$-\infty$.
Therefore not only is $c_a^{\ast} > 0$, at $c=c_a^{\ast}$ the minimiser of
(\ref{corr2}) is discontinuous. Heuristically, this means that the optimal
strategy
stays localised on scale $t^{1/d}$ as $c\da c_a^{\ast}$, i.e., the collapse
transition is {\it first order}.

The analogue of Corollary \ref{c1} for $d=2$ follows from Theorem \ref{t2}
and reads as follows.

\vglue9pt

\specialnumber{2} \proclaim{{C}orollary} \label{c2}
Let $d=2$ and $a>0$. For every $c>0,$
\be{corr5}
\lim\limits_{t \ra \infty} \frac{1}{\log t}
\log E\biggl(\exp[-ct^{-1}\log^2 t~|W^a(t)|]\biggr) = -J^{2\pi}(c)
\ee
with
\be{corr6}
J^{2\pi}(c) = \inf\limits_{b>0} [bc + I^{2\pi}(b)].
\ee
\ec
The same statements as in (\ref{corr3}), (\ref{corr4}) hold, again with
$c^{\ast}>0$,
because by Theorem 4(ii) the left derivative of $I^{2\pi}$ at $b=2\pi$ is
strictly negative.
Thus, also for $d=2$ there is a collapse transition. However, $b \mapsto
I^{2\pi}(b)/(2\pi-b)$
is strictly decreasing on $(0,2\pi)$ by Theorems \ref{t3}(i) and
\ref{t4}(ii), and so
at $c=c^{\ast}$ the minimiser is continuous. This means that the optimal
strategy does
not stay localised on scale $\sqrt{t/\log t}$ as $c \da c^{\ast}$, i.e.,
the collapse
transition is {\it second order}.

\demo{{\rm 1.7.} Upward deviations} Finally, the moderate and large deviations of
$|W^a(t)|$ in the {\it
upward} direction are a complicated issue. Here the optimal strategy is
entirely different
from the previous ones, because the Wiener sausage tries to expand rather
than to
contract. Partial results have been obtained by van den Berg and T\'oth
\cite{BeT} and
van den Berg and Bolthausen \cite{BeBo}.\enddemo

\section{The upper bound in Theorem \ref{t1}}

This section contains the main probabilistic part of the paper and,
together with Sections 3 and 4,
provides the proof of Theorems \ref{t1} and \ref{t2}.

\demo{{\rm 2.1.} Compactification\/{\rm :} Propositions {\rm \ref{p1}} and {\rm \ref{p2}}}
We begin by doing a standard compactification. Let $\Lambda_N$ be the torus
of size $N>0$, i.e., $[-\frac{N}{2},\frac{N}{2})^d$ with periodic boundary
conditions.
For $t>0$, let
$\beta_{Nt^{1/d}}(s)$, $s \geq 0$, be the Brownian motion wrapped around
$\Lambda_{Nt^{1/d}}$, and let $W^a_{Nt^{1/d}}(s)$, $s \geq 0$, denote its
Wiener
sausage. Then trivially
\be{ubper}
P(|W^a(t)| \leq bt) \leq P(|W^a_{Nt^{1/d}}(t)| \leq bt)
\ee
for all $a>0$, $b>0$, $N>0$ and $t>0$. Next, by Brownian scaling,
$|W^a_{Nt^{1/d}}(t)|$
has the same distribution as $t|W^{at^{-1/d}}_N(t^{(d-2)/d})|$. Hence, putting
\be{tau}
\tau = t^{(d-2)/d}
\ee
we get
\be{ubscal}
P(|W^a(t)| \leq bt) \leq P(|W^{a\tau^{-1/(d-2)}}_N(\tau)| \leq b).
\ee
The right-hand side of (\ref{ubscal}) involves the Wiener sausage on
$\Lambda_N$ at time $\tau$ with a radius that shrinks with $\tau$.

In Sections 2.2--2.5 we shall prove the following: 

\specialnumber{1} \proclaim{Proposition} \label{p1}
Let $d \geq 3$ and $a>0${\rm .} For every $b>0$ and $N>0,$
\be{mdsshrink}
\lim_{\tau \ra \infty} \frac{1}{\tau} \log P(|W^{a\tau^{-1/(d-2)}}_N(\tau)|
\leq b) = - I^{\kappa_a}_N(b),
\ee
where $I^{\kappa_a}_N(b)$ is given by the same formulas as in
{\rm (\ref{mdsrf}), (\ref{mdsset}),}
except that $\R^d$ is replaced by $\Lambda_N${\rm .}
\ep

{}From (\ref{tau})--(\ref{mdsshrink}) we get
\be{mdsunshrink}
\limsup_{\tau \ra \infty} \frac{1}{t^{(d-2)/d}} \log P(|W^a(t)| \leq bt)
\leq -I^{\kappa_a}_N(b)
\qquad \mbox{ for all } N>0.\hskip.25in
\ee
In Section 2.6 we shall show:
\pagebreak
\specialnumber{2} \proclaim{Proposition} \label{p2}
$\lim\limits_{N \ra \infty} I^{\kappa_a}_N(b) = I^{\kappa_a}(b)$ for all
$a>0$ and $b>0${\rm .}
\ep

\noindent Combining this with (\ref{mdsunshrink}) we get the upper bound in Theorem
\ref{t1}.

Our proof of Proposition \ref{p1} is based on a new approach for treating
large deviations
of the Wiener sausage on the torus. This approach uses a conditioning
argument, a version of Talagrand's concentration inequality, and the most
basic LDP  (Large Deviation Principle) of Donsker and Varadhan.

Throughout the rest of this section the Brownian motion lives on
$\Lambda_N$ with
$N$ fixed, and we suppress the indices $a$ and $N$ from most expressions.
Abbreviate
\be{Vdef}
V_{\tau}= |W^{a\tau^{-1/(d-2)}}_N(\tau)|.
\ee
We shall prove the following:

\specialnumber{3} \proclaim{Proposition} \label{VLDP}
$(V_{\tau})_{\tau>0}$ satisfies the {\rm LDP} on $\R_+$ with rate $\tau$ and with
rate
function
\be{J}
J^{\kappa_a}_N(b)=\inf_{\phi \in \partial \Phi^{\kappa_a}_N(b)}\left[
\frac{1}{2} \int_{\Lambda_N} |\nabla\phi|^2(x) dx \right]
\ee
with
\be{Jset}
\partial \Phi^{\kappa_a}_N(b) = \left\{\phi\in H^1(\Lambda_N) \colon \,
\int_{\Lambda_N} \phi^2(x)dx = 1, \,\int_{\Lambda_N}
\Big(1-e^{-\kappa_a \phi^2(x)}\Big)dx = b \right\}.
\ee
\ep

Proposition \ref{VLDP} obviously implies Proposition \ref{p1}. We shall
see in Section 3 that it is also the key to the lower bound in Theorem
\ref{t1}, but this requires a separate argument.

The form of Proposition \ref{VLDP} suggests that some kind of contraction
principle is in force. However, it seems to be impossible to approach the
problem directly from that angle. Instead, we use an approximation argument
consisting of three steps:
 
\smallbreak\noindent \hangindent=28pt\hangafter=1 \hskip16pt $\bullet$ \hskip3pt  Section 2.2:
For $\varepsilon>0$, 

\phantom{time to come home}
\vglue-24pt
\be{barXdef}
\X_{\tau,\varepsilon} = \{\beta(i\varepsilon)\}_{1 \leq i \leq \tau/\varepsilon}.
\ee
\noindent \hangindent=28pt\hangafter=0  (For notational convenience $\tau/\varepsilon$ is taken to be integral.)
We first approximate $V_{\tau}$ by $\E_{\tau,\varepsilon}(V_{\tau})$,
where $\E_{\tau,\varepsilon}$ denotes the {\it conditional\/} expectation
given $\X_{\tau,\varepsilon}$. We prove that the difference between $V_{\tau}$
and $\E_{\tau,\varepsilon}(V_{\tau})$ is negligible in the limit as $\tau \to
\infty$
followed by $\varepsilon \da 0$. This is done by   application of a concentration
inequality of Talagrand.

\smallbreak\noindent \hangindent=28pt\hangafter=1 \hskip16pt $\bullet$ \hskip3pt Section 2.3:
We represent $\E_{\tau,\varepsilon}(V_{\tau})$ as a functional of the empirical
measure

\phantom{time to come home}
\vglue-30pt
\be{empdef}
L_{\tau,\varepsilon} = \frac{\varepsilon}{\tau} \sum_{i=1}^{\tau/\varepsilon}
\delta_{\bigl(\beta((i-1)\varepsilon),\beta(i\varepsilon)\bigr)}.
\ee
\noindent \hangindent=28pt\hangafter=0  According to Donsker and Varadhan, $(L_{\tau,\varepsilon})_{\tau>0}$
satisfies an LDP. We need some further approximations to get the dependence
of $\E_{\tau,\varepsilon}(V_{\tau})$ on $L_{\tau,\varepsilon}$ in a suitable
form, but essentially based on just this LDP we get an LDP for
$(\E_{\tau,\varepsilon}(V_{\tau}))_{\tau>0}$ via a contraction principle.

\medbreak\noindent \hangindent=28pt\hangafter=1 \hskip16pt $\bullet$ \hskip3pt Section 2.4:
Finally, we have to perform the limit $\varepsilon\da 0$. By our previous
result we already know that $V_{\tau}$ is well approximated by
$\E_{\tau,\varepsilon}(V_{\tau})$. It therefore suffices to have
an appropriate approximation for the variational formula in the LDP
for $(\E_{\tau,\varepsilon}(V_{\tau}))_{\tau>0}$.
\bigbreak

\noindent In Section 2.5 the above results are collected to prove Proposition
\ref{VLDP}. 

It will be expedient to use the abbreviation
\be{Tdef}
T_{\tau} = \tau^{2/(d-2)}.
\ee
So the radius of our Wiener sausage on $\Lambda_N$ is $a/\sqrt{T_{\tau}}$.

\vglue6pt

\demo{{\rm 2.2.} Approximation of $V_{\tau}$ by $\E_{\tau,\varepsilon}(V_{\tau})$} 
Recall the definition of $\X_{\tau,\varepsilon}$ in (\ref{barXdef}).
We denote by $\P_{\tau,\varepsilon}$ and $\E_{\tau,\varepsilon}$ the
conditional probability and expectation given $\X_{\tau,\varepsilon}$.

The main result of this section is that $V_{\tau}$ is well approximated
by $\E_{\tau,\varepsilon}(V_{\tau})$ in the following sense: \enddemo

\specialnumber{4} \proclaim{Proposition} \label{Concentration}
For all $\delta>0,$
\be{C1}
\lim_{\varepsilon \da 0} \limsup_{\tau\rightarrow\infty}
\frac{1}{\tau}\log P(|V_{\tau}-\E_{\tau,\varepsilon}(V_{\tau})|
\geq \delta) = -\infty.
\ee
\ep

\demo{Proof}
The proof proceeds via a series of estimates.
\smallbreak
{1.} We begin by truncating the excursions. Define
\be{Wurstscheiben}
W_i = \bigcup_{(i-1)\varepsilon \leq s \leq i\varepsilon}
B_{a/\sqrt{T_{\tau}}}(\beta(s))
\qquad (1 \leq i \leq \tau/\varepsilon).
\ee
Then
\be{Vdef*}
V_{\tau}=\left|\bigcup\limits_{i=1}^{\tau/\varepsilon}W_i\right|.
\ee
For $K>0$, let
\be{exctr}
J^K_{\tau,\varepsilon} = \{1 \leq i \leq \tau/\varepsilon \colon
~|\beta((i-1)\varepsilon)-\beta(i\varepsilon)| \leq K\sqrt{\varepsilon}\,\}
\ee
and define
\be{VKdef}
V^K_{\tau,\varepsilon} = \left|\bigcup\limits_{i=1}^{\tau/\varepsilon}
W_i 1\{i \in J^K_{\tau,\varepsilon}\}\right|, \qquad
\widehat V^K_{\tau,\varepsilon} = \left|\bigcup\limits_{i=1}^{\tau/\varepsilon}
W_i 1\{i \notin J^K_{\tau,\varepsilon}\}\right|.\hskip.5in
\ee
Since $0 \leq V_\tau - V^K_{\tau,\varepsilon} \leq \widehat
V^K_{\tau,\varepsilon}$, we have
\be{trin}
|V_{\tau} - \E_{\tau,\varepsilon}(V_{\tau})| \leq
|V^K_{\tau,\varepsilon} - \E_{\tau,\varepsilon}(V^K_{\tau,\varepsilon})|
+ \widehat V^K_{\tau,\varepsilon} +\E_{\tau,\varepsilon}(\widehat
V^K_{\tau,\varepsilon}).\hskip.5in
\ee
The claim will follow after we prove the following two results:
\begin{eqnarray} \label{alrlim}
&&\\
\lim\limits_{\varepsilon \da 0} \limsup\limits_{\tau\rightarrow\infty}
\frac{1}{\tau}\log
P(|V^K_{\tau,\varepsilon}-\E_{\tau,\varepsilon}(V^K_{\tau,\varepsilon})|
\geq \delta) &\hskip-4pt = \hskip-4pt & -\infty \enspace \mbox{ for all } \delta>0, K \geq
K_0(\delta)\nonumber \\[4pt]
\lim\limits_{\varepsilon \da 0} \limsup\limits_{\tau\rightarrow\infty}
\frac{1}{\tau}\log P(\widehat V^K_{\tau,\varepsilon}\geq\delta) &\hskip-4pt = \hskip-4pt & -\infty \enspace
\mbox{ for all } \delta>0, K \geq K_0(\delta).\nonumber
\end{eqnarray}
Indeed, the third term on the right-hand side of (\ref{trin}) needs no
extra consideration,
because $\widehat V^K_{\tau,\varepsilon} \leq |\Lambda_N|$ implies that
$\E_{\tau,\varepsilon}
(\widehat V^K_{\tau,\varepsilon}) \leq \frac{\delta}{2} +
|\Lambda_N|\P_{\tau,\varepsilon}
(\widehat V^K_{\tau,\varepsilon} \geq \frac{\delta}{2})$ and hence
\be{alrlim*}
P\Big(\E_{\tau,\varepsilon}(\widehat V^K_{\tau,\varepsilon}) \geq \delta\Big)
\leq P\Big(\P_{\tau,\varepsilon}(\widehat V^K_{\tau,\varepsilon} \geq
\textstyle{\frac{\delta}{2}})
\geq \textstyle{\frac{\delta}{2|\Lambda_N|}}\Big) \leq
\textstyle{\frac{2|\Lambda_N|}{\delta}}
P\Big(\widehat V^K_{\tau,\varepsilon} \geq \textstyle{\frac{\delta}{2}}\Big).
\ee
\bigbreak
2.  To prove the second claim in (\ref{alrlim}), we estimate
\begin{eqnarray} \label {Vhat1} &&\\
\noalign{\vskip-8pt}
  P(\widehat V^K_{\tau,\varepsilon}\geq\delta)
&\leq& e^{-\delta\tau/2\varepsilon} E\Big(\exp\left[\frac{\tau}{2\varepsilon}
\sum\limits_{i=1}^{\tau/\varepsilon} |W_i|1\{i \notin
J^K_{\tau,\varepsilon}\}\right]\Big)\nonumber\\  
&=& e^{-\delta\tau/2\varepsilon} \Big\{E\Big(\exp\left[\frac{\tau}{2\varepsilon}
|W_1|1\{1 \notin J^K_{\tau,\varepsilon}\}\right]\Big)\Big\}^{\tau/\varepsilon}\nonumber
\\
&=& e^{-\delta\tau/2\varepsilon} \Big\{1+E\Big(1\{1 \notin J^K_{\tau,\varepsilon}\}
\Big(\exp\left[\frac{\tau}
{2\varepsilon}|W_1|\right]\Big)-1\Big)\Big\}^{\tau/\varepsilon}\nonumber \\ 
&\leq& e^{-\delta\tau/2\varepsilon} \Big\{1+\sqrt{\delta_K
C_{\tau,\varepsilon}}\Big\}^{\tau/\varepsilon},\nonumber
\end{eqnarray}
where
\begin{eqnarray} \label{Vhat2}
\delta_K &=& P(|\beta(\varepsilon)|>K\sqrt{\varepsilon})\\[0.4cm]
C_{\tau,\varepsilon} &=& E\Big(\exp\left[\frac{\tau}{\varepsilon}|W^{a/\sqrt{T_\tau}}
(\varepsilon)|\right]\Big).\nonumber
\end{eqnarray}
It is evident that $|W^{a/\sqrt{T_\tau}}(\varepsilon)|$ is smaller on the
torus than on $\R^d$.
Therefore we get, after Brownian scaling and  using that
$\tau/T_\tau^{d/2}=1/T_\tau$ by (\ref{Tdef}),
\be{Vhat3}
C_{\tau,\varepsilon} \leq E\Big(\exp\left[\frac{1}{\varepsilon T_\tau}
\left|W^a(\varepsilon T_\tau)\right|\right]\Big).
\ee
It follows from the results in van den Berg and Bolthausen \cite{BeBo} that
\be{boundheat}
\sup\limits_{T \geq 1} E\Big(\exp\left[\frac{1}{T}|W^a(T)|\right]\Big) <
\infty.
\ee
Hence the right-hand side of (\ref{Vhat3}) is bounded above by some $C<\infty$
for all $\tau \geq \tau_0(\varepsilon)$, and so we find that
\be{Vhat4}
\limsup_{\tau\rightarrow\infty}
\frac{1}{\tau}\log P(\widehat V^K_{\tau,\varepsilon} \geq \delta) \leq
- \frac{\delta}{2\varepsilon} + \frac{\sqrt{\delta_K C}}{\varepsilon}
\qquad \mbox{ for all } \varepsilon,K>0.\hskip.25in
\ee
Since $\lim_{K \to \infty} \delta_K=0$, there exists a $K_0(\delta)$ such that
$\sqrt{\delta_K C} \leq \frac{\delta}{4}$ for $K \geq K_0(\delta)$. For such
$K$ we now let $\varepsilon \da 0$ to get the second claim in (\ref{alrlim}).\bigbreak

 3.  To prove the first claim in (\ref{alrlim}) we argue as follows.
Conditionally on $\X_{\tau,\varepsilon}$, the $W_{i}$ are independent random
open subsets of $\Lambda_{N}$. Let $S$ be the set of open subsets of $\Lambda_{N}$.
The mapping $d \colon S \times S \mapsto [0,\infty)$ with $d(A,B)=\left|A
\Delta B\right|$
defines a pseudometric on $S$. We equip $S$ with the Borel field ${\frak S}$
generated by this pseudometric. Then $\P_{\tau,\varepsilon}$ defines a product
measure on $(S,{\frak S})^{\tau/\varepsilon}$, which we denote by the same
symbol $\P_{\tau,\varepsilon}$. Define
\be{1}
V(C) = \left|\bigcup\limits_{i \in J^K_{\tau,\varepsilon}} C_{i}\right|
\qquad (C=\{C_{i}\} \in S^{\tau/\varepsilon})
\ee
(note that $\X_{\tau,\varepsilon}$ fixes $J^K_{\tau,\varepsilon}$). Clearly, $V$
is Lipschitz
in the sense that
\be{2}
|V(C)-V(C')|  \leq \sum_{i \in J^K_{\tau,\varepsilon}}|C_{i}\Delta C'_{i}|
\qquad (C,C' \in S^{\tau/\varepsilon}).
\ee

\bigbreak 4.  Let us denote by $m^K_{\tau,\varepsilon}$ the median of the
distribution of
$V^K_{\tau,\varepsilon}$ under the conditional law $\P_{\tau,\varepsilon}$. Define
\be{3}
A=\{C\in S^{\tau/\varepsilon} \colon ~V(C)\leq m^K_{\tau,\varepsilon}\}.
\ee
Since the distribution of $V^K_{\tau}$ under $\P_{\tau,\varepsilon}$ has no
atoms, we
have $\P_{\tau,\varepsilon}(A)=\frac{1}{2}$. From Talagrand \cite[Th.\ 2.4.1]{T}
(see also Remark 2.1.3) we therefore have
\be{4}
\E_{\tau,\varepsilon}\Big(\exp[\lambda f(A,\{W_{i}\})]\Big)
\leq 2 \prod_{i \in J^K_{\tau,\varepsilon}}\E_{\tau,\varepsilon}
\Big(\cosh[\lambda|W_{i}\Delta W_{i}^{\prime}|]\Big),\hskip.5in
\ee
where
\be{fdef}
f(A,\{C_{i}\})=\inf_{C'\in A} \sum_{i\in J^K_{\tau,\varepsilon}}|C_{i}\Delta
C'_{i}|
\ee
and $\{W_{i}^{\prime}\}$ is an independent copy of $\{W_{i}\}$. From the Markov
inequality we therefore get  
\be{Phi**}
\P_{\tau,\varepsilon}(f(A,\{W_{i}\}) \geq \delta) \leq
2 \inf_{\lambda>0}e^{-\lambda\delta}\prod_{i \in J^K_{\tau,\varepsilon}}
\E_{\tau,\varepsilon}\Big(\cosh[\lambda|W_{i}\Delta W_{i}^{\prime}|]\Big)
=:\Xi^K_{\tau,\varepsilon}(\delta).
\end{equation}
Arguing similarly with $\widehat A=\left\{C\in S^{\tau/\varepsilon}\colon
~V(C)\geq m^K_{\tau,\varepsilon}\right\}$ we get, by \hbox{(\ref{2}), (\ref{3}),}
\be{5}
\P_{\tau,\varepsilon}(\left|V^K_{\tau,\varepsilon}-m^K_{\tau,\varepsilon}\right|
\geq \delta)
\leq \P_{\tau,\varepsilon}(f(A,\{W_{i}\}) \geq \delta)
+ \P_{\tau,\varepsilon}(f(\widehat A,\{W_{i}\}) \geq \delta)
\leq 2\Xi^K_{\tau,\varepsilon}(\delta).
\ee

\medbreak 5.  Next, since $V^K_{\tau,\varepsilon} \leq \left|\Lambda_{N}\right|$ we have
\be{6}
\left|\E_{\tau,\varepsilon}(V^K_{\tau,\varepsilon})-m^K_{\tau,\varepsilon}\right| \leq
\textstyle{\frac{\delta}{3}}+\left|\Lambda_{N}\right|\P_{\tau,\varepsilon}
\left(\left|V^K_{\tau,\varepsilon}-m^K_{\tau,\varepsilon}\right| \geq
\textstyle{\frac{\delta}{3}}\right)
\ee
and consequently
\begin{eqnarray} \label{7}
\quad&& \P_{\tau,\varepsilon}\left(\left|V^K_{\tau,\varepsilon}-\E_{\tau,\varepsilon}
(V^K_{\tau,\varepsilon})\right|\geq\delta\right)\\
&&\hskip.5in \leq
\P_{\tau,\varepsilon}\left(\left|V^K_{\tau,\varepsilon}-m^K_{\tau,\varepsilon}\right|
\geq \frac{\delta}{3}\right) +
1\left\{\left|\E_{\tau,\varepsilon}(V^K_{\tau,\varepsilon})
-m^K_{\tau,\varepsilon}\right|\geq \frac{2\delta}{3}\right\}\nonumber  \\
&&\hskip.5in \leq 2\Xi^K_{\tau,\varepsilon}\left(\frac{\delta}{3}\right)
+
1\left\{\P_{\tau,\varepsilon}\left(\left|V^K_{\tau,\varepsilon}-m^K_{\tau,\varepsilon
}\right|
\geq\frac{\delta}{3}\right)
\geq\frac{\delta}{3\left|\Lambda_{N}\right|}\right\} \nonumber  \\
&&\hskip.5in\leq 2\Xi^K_{\tau,\varepsilon}\left(\frac{\delta}{3}\right)
+ 1\left\{2\Xi^K_{\tau,\varepsilon}\left(\frac{\delta}{3}\right)
\geq\frac{\delta}{3\left|\Lambda_{N}\right|}\right\}.
\nonumber
\end{eqnarray}
Using this inequality we get, after averaging over $\X_{\tau,\varepsilon}$,
\be{8}
P(\left|V^K_{\tau,\varepsilon}-\E_{\tau,\varepsilon}(V^K_{\tau})\right| \geq \delta)
\leq \left(1+\textstyle{\frac{3\left|\Lambda_{N}\right|}{\delta}}\right)
E\left(2\Xi^K_{\tau,\varepsilon}\left(\textstyle{\frac{\delta}{3}}\right)\right).
\ee
In order to prove the first claim in (\ref{alrlim}), it therefore suffices
to show that
\be{basic}
\lim_{\varepsilon \da 0}\limsup_{\tau\rightarrow\infty}
\frac{1}{\tau}\log E\left(\Xi^K_{\tau,\varepsilon}(\delta)\right)
= -\infty \qquad \mbox{ for all } \delta>0,K \geq
K_0(\delta). \hskip.5in
\ee
We shall actually prove more, namely that
\be{basic2}
\lim_{\varepsilon \da 0}\limsup_{\tau\rightarrow\infty}
\frac{1}{\tau}\log\left\|\Xi^K_{\tau,\varepsilon}(\delta)\right\|_{\infty}
=-\infty \qquad \mbox{ for all } \delta>0,K \geq
K_0(\delta).\hskip.5in
\ee

\bigbreak 6.  In order to estimate $\E_{\tau,\varepsilon}(\cosh[\lambda\left|W_{i}
\Delta W_{i}^{\prime}\right|])$ in (\ref{Phi**}), we pick
$\lambda=c\tau/\varepsilon$
with $0<c\leq 1$ and use the fact that $\cosh(cd)\leq1+c^{2}\exp(d)$ for $0<c\leq1$
and $d>0$. For $x\in\Lambda_{N}$, we write $E_{x,\varepsilon}$ to denote
expectation under a Brownian bridge of length $\varepsilon$ between $0$ and $x$,
i.e., a Brownian motion starting at $0$ and conditioned to be at $x$ at
time $\varepsilon$.
It is evident that the volume of the Wiener sausage associated with such a
Brownian
bridge is smaller on the torus than on $\R^d$. Thus we have
\begin{eqnarray} \label{9}
\quad &&\E_{\tau,\varepsilon}\Big(\cosh\left[c\frac{\tau}{\varepsilon}
\left|W_{i}\Delta W_{i}^{\prime}\right|\right]\Big)\\[4pt]
&&\hskip1in \leq 1+c^2\E_{\tau,\varepsilon}\Big(\exp\left[\frac{\tau}{\varepsilon}
\left|W_{i}\Delta W_{i}^{\prime}\right|\right]\Big)\nonumber \\[4pt]
&&\hskip1in\leq  1+c^2\left\{\E_{\tau,\varepsilon}\Big(\exp\left[\frac{\tau}{\varepsilon}
\left|W_i\right|\right]\Big)\right\}^2\nonumber\\[4pt]
&&\hskip1in\leq  1+c^2\left\{\sup\limits_{|x|\leq K}E_{x\sqrt{\varepsilon},\varepsilon}
\left(\exp\left[\frac{\tau}{\varepsilon}\left|W^{a/\sqrt{T_{\tau}}}
(\varepsilon)\right|\right]\right)\right\}^{2},\nonumber
\end{eqnarray}
where we recall (\ref{Wurstscheiben}) and use the fact that
$|\beta((i-1)\varepsilon)-\beta(i\varepsilon)|
\leq K\sqrt{\varepsilon}$ for $i \in J^K_{\tau,\varepsilon}$. By Brownian scaling
we get  
($\tau/T_{\tau}^{d/2}=1/T_{\tau}$) \vglue-9pt
\be{10}
E_{x\sqrt{\varepsilon},\varepsilon}\left(\exp\left[\frac{\tau}{\varepsilon}
\left|W^{a/\sqrt{T_{\tau}}}(\varepsilon)\right|\right]\right)\\ 
\leq E_{x\sqrt{\varepsilon T_{\tau}},\varepsilon T_{\tau}
}\left(\exp\left[\frac{1}{\varepsilon T_{\tau} }
\left|W^{a}(\varepsilon T_{\tau} )\right|\right]\right).
\ee

\bigbreak 7. With the help of Lemma \ref{exp_moment} below it follows from
(\ref{10}) that there
exists a $C_K<\infty$ such that for all $\tau\geq\tau_0(\varepsilon),$
\be{11}
\sup\limits_{|x|\leq
K}E_{x\sqrt{\varepsilon},\varepsilon}\left(\exp\left[\frac{\tau}{\varepsilon}
\left|W^{a/\sqrt{T_{\tau}}}(\varepsilon)\right|\right]\right)
\leq C_K.
\ee
Therefore, combining (\ref{Phi**}), (\ref{9}) and (\ref{11}) we get
\begin{eqnarray} \label{N1}
\Xi^K_{\tau,\varepsilon}(\delta) &\leq &2e^{-c\delta\frac{\tau}{\varepsilon}}
\prod_{i \in J^K_{\tau,\varepsilon}}\E_{\tau,\varepsilon}
\Big(\cosh\left[c\frac{\tau}{\varepsilon}\left|W_{i}\Delta
W_{i}^{\prime}\right|\right]\Big)\\
&\leq& 2e^{-c\delta\frac{\tau}{\varepsilon}}
\prod_{i=1}^{\tau/\varepsilon}(1+c^{2}C_K^{2})
\leq 2e^{(-c\delta+c^{2}C_K^{2})\frac{\tau}{\varepsilon}}.\nonumber
\end{eqnarray}
Pick $c=\delta/2C_K^2$ and note that there exists a $K_0(\delta)$ such that
$0<c\leq 1$ for
$K\geq K_0(\delta)$. Let $\tau\ra\infty$ followed by $\varepsilon\da 0$, to
get (\ref{basic2}).
The proof of Proposition \ref{Concentration} is now complete.
\enddemo 

We conclude this section with the following fact:
\specialnumber{1} \proclaim{Lemma} \label{exp_moment}
For every $K>0$ there exists a $C_K<\infty$ such that
\be{WSexpest}
\sup_{T\geq 2} \sup_{\left|x\right|\leq K}
E_{x\sqrt{T},T}\left(\exp\left[\frac{1}{T}\left|W^{a}(T)\right|\right]\right)
\leq C_K.
\ee
\el

\demo{Proof}
We begin by removing the bridge restriction. Write $$p_t(x,y) = (2\pi t)^{-d/2} 
 \exp[-|x-y|^2/2t]$$ to denote the heat kernel on $\R^d$ and put
$p_t(x)=p_t(0,x)$.
Write $E_{y,t;z,2t}$ to denote expectation under a Brownian motion starting
at 0
and conditioned to be at $y$ at time $t$ and at $z$ at time $2t$. Then we
may estimate
\begin{eqnarray}\label{heatest*}
&&E_{z,2t}\Big(\exp\Big[\frac{1}{2t}|W^a(2t)|\Big]\Big)\\ [0.6pt]
&&\hskip.75in = \int_{\R^d} dy ~\frac{p_t(y)p_t(z-y)}{p_{2t}(z)}
~E_{y,t;z,2t}\Big(\exp\Big[\frac{1}{2t}|W^a(2t)|\Big]\Big)  \nonumber \\ [0.6pt]
&&\hskip.75in \leq \frac{1}{p_{2t}(z)} \int_{\R^d} dy
~p_t(y)~E_{y,t}\Big(\exp\Big[\frac{1}{2t}|W^a(t)|\Big]\Big)\nonumber \\[0.6pt]
&&\hskip.75in  \quad \times
~p_t(z-y)~E_{z-y,t}\Big(\exp\Big[\frac{1}{2t}|W^a(t)|\Big]\Big)\nonumber  \\[0.6pt]
&&\hskip.75in \leq \frac{1}{p_{2t}(z)} \int_{\R^d} dy ~\Big\{p_t(y)
~E_{y,t}\Big(\exp\Big[\frac{1}{2t}|W^a(t)|\Big]\Big)\Big\}^2\nonumber  \\[0.6pt]
&&\hskip.75in \leq \frac{p_t(0)}{p_{2t}(z)} \int_{\R^d} dy ~p_t(y)
~E_{y,t}\Big(\exp\Big[\frac{1}{t}|W^a(t)|\Big]\Big)\nonumber  \\[0.6pt]
&&\hskip.75in = \frac{p_t(0)}{p_{2t}(z)} ~E\Big(\exp\Big[\frac{1}{t}|W^a(t)|\Big] \Big).
\nonumber 
\end{eqnarray}
Here we use, respectively, the subadditivity of $t \mapsto |W^a(t)|$,
Cauchy-Schwarz,
Jensen and the bound $p_t(y)\leq p_t(0)$. Next put $z=x\sqrt{T}$, $t=T/2$
in (\ref{heatest*})
and use $\sup_{|x|\leq K} p_{T/2}(0)/p_{T}(x\sqrt{T})$
$=2^{d/2}\exp(K^2/2)$, to see that
the claim follows from (\ref{boundheat}). $\phantom{\sum^|}$
\enddemo

\demo{{\rm 2.3.} The {\rm LDP} for $(\E_{\tau,\varepsilon}(V_{\tau}))_{\tau>0}$}
Let $I_{\varepsilon}^{(2)} \colon {\cal M}^+_1(\Lambda_N \times \Lambda_N)
\mapsto [0,\infty]$ be the entropy function
\be{27}
I_{\varepsilon}^{(2)}(\mu)=\left\{
\begin{array}
[c]{cc}
h(\mu|\mu_{1}\otimes\pi_{\varepsilon}) &\mbox{if }\mu_{1}=\mu_{2}\\
\infty &\mbox{otherwise},
\end{array}
\right.
\ee
where $h(\cdot|\cdot)$ denotes relative entropy between measures, $\mu_{1}$
and $\mu_{2}$
are the two marginals of $\mu$, and $\pi_{\varepsilon}(x,dy)=p_{\varepsilon}(y-x)dy$
is the Brownian transition kernel on $\Lambda_N$. Furthermore, for $\eta>0$
let $\Phi_{\eta}
\colon {\cal M}^+_1(\Lambda_N \times \Lambda_N) \mapsto [0,\infty)$ be
the function
\be{Phiext}
\Phi_{\eta}(\mu)=\int_{\Lambda_{N}}dx
\left(1-\exp\left[-\eta\kappa_{a}\int_{\Lambda_N\times\Lambda_N}
\varphi_{\varepsilon}(y-x,z-x)
\mu(dy,dz)\right]\right) \enspace
\ee
with
\be{15}
\varphi_{\varepsilon}(y,z)=\frac{\int_{0}^{\varepsilon}
ds\,p_{s}(-y)p_{\varepsilon-s}(z)}
{p_{\varepsilon}(z-y)}.
\ee
Our main result in this section is the following:

\specialnumber{5} \proclaim{Proposition} \label{Epsilon_LDP}
$(\E_{\tau,\varepsilon}(V_{\tau}))_{\tau>0}$ satisfies the {\rm LDP} on $\R_+$ with
rate $\tau$ and
with rate function
\be{28}
J_{\varepsilon}(b)=\inf\left\{\frac{1}{\varepsilon}I_{\varepsilon}^{(2)}(\mu)
\colon \,\mu \in {\cal M}^+_1(\Lambda_N \times \Lambda_N), \,
\Phi_{1/\varepsilon}(\mu)=b\right\}.\hskip.5in
\ee
\ep

\demo{{P}roof}
Throughout the proof, $c_1,c_2,\dots$ are constants that may depend on
$a,\varepsilon,N$
(which are fixed) but not on any of the other variables.\medbreak
{1.} First we approximate $V_{\tau}$ by cutting out small holes around
the points $\beta(i\varepsilon), 1\leq i\leq\tau/\varepsilon$. Fix $K>0$, let
\be{18}
W_{i}^{K/\sqrt{T_{\tau}}}=W_{i} \Big\backslash \left[B_{K/\sqrt{T_{\tau}}}
\Big(\beta((i-1)\varepsilon)\Big)\cup
B_{K/\sqrt{T_{\tau}}}\Big(\beta(i\varepsilon)\Big)\right]\hskip.5in
\ee
and put
\be{19}
V_{\tau}^{K}=\left|\bigcup_{i=1}^{\tau/\varepsilon}
W_{i}^{K/\sqrt{T_{\tau}}}\right|.
\ee
Clearly, we have cut out at most $\tau/\varepsilon+1$ times the volume of a
ball of radius $K/\sqrt{T_{\tau}}$, so
\be{cut}
\left|V_{\tau}-V_{\tau}^{K}\right|\leq c_1 K^{d}/T_{\tau},
\end{equation}
which tends to zero as $\tau\ra\infty$ and therefore is negligible for our
purpose. This
cutting procedure is convenient as will become clear later on.\medbreak

{2.} For $y,z\in\Lambda_{N}$, define
\be{14}
q_{\tau,\varepsilon}(y,z) = P_{y,z}(\sigma_{a/\sqrt{T_\tau}} \leq \varepsilon),
\ee
where $P_{y,z}(\cdot) = P((\beta(t))_{t\in[0,\varepsilon]} \in \cdot \mid
\beta(0)=y,
\beta(\varepsilon)=z)$ and $\sigma_{a/\sqrt{T_\tau}}$ is the first entrance
time into
$B_{a/\sqrt{T_\tau}}=B_{a/\sqrt{T_\tau}}(0)$. We can express
$\E_{\tau,\varepsilon}(V_{\tau})$
in terms of $q_{\tau,\varepsilon}(y,z)$ and the empirical measure
$L_{\tau,\varepsilon}$ defined
in (\ref{empdef}) as follows:
\bel{repr}
 && \\
 \E_{\tau,\varepsilon}(V_{\tau}^{K})&\hskip-12pt =\hskip-12pt&
\int_{\Lambda_{N}}dx\left(1-\P_{\tau,\varepsilon}\left(x\notin
\bigcup_{i=1}^{\tau/\varepsilon}W_{i}^{K/\sqrt{T_{\tau}}}\right)\right)\nn  \\
&\hskip-12pt =\hskip-12pt&\int_{\Lambda_{N}}dx\left(1-\prod\limits_{i=1}^{\tau/\varepsilon}
\left\{1-\P_{\tau,\varepsilon}\left(x\in
W_{i}^{K/\sqrt{T_{\tau}}}\right)\right\}\right)\nn \\
&\hskip-12pt =\hskip-12pt&\int_{\Lambda_{N}}dx\left(1-\exp\left[\frac{\tau}{\varepsilon}
\int_{\Lambda_N\times\Lambda_N}\right.\right.\nn\\
&&\hskip.92in \left.\left.\phantom{\int^1} \log\!
\left(1-q_{\tau,\varepsilon}^{K/\sqrt{T_{\tau}}}
(y-x,z-x)\right)
L_{\tau,\varepsilon}(dy,dz)\right]\right)\!,
\nne
where for $\rho>0$ we define
$q_{\tau,\varepsilon}^{\rho}(y,z)=q_{\tau,\varepsilon}(y,z)$ if
$y,z\notin B_{\rho}$ and zero otherwise.\medbreak

3.  We want to expand the logarithm and do an approximation. For this
we need the following facts about Brownian motion on $\Lambda_N$, which
come as an intermezzo.
Recall that $\kappa_{a}$ is the Newtonian capacity
of the ball with radius~$a$.
\enddemo

\specialnumber{2} \proclaim{Lemma} \label{cut1}
{\rm (a)} $\lim_{K\rightarrow\infty}\limsup_{\tau\rightarrow\infty}
\sup_{y,z\notin B_{K/\sqrt{T_{\tau}}}}q_{\tau,\varepsilon}(y,z)=0${\rm .}
\medbreak
{\rm (b)} $\lim_{\tau\rightarrow\infty}\sup_{y,z\notin B_{\rho}}
\left|\tau q_{\tau,\varepsilon}(y,z)-\kappa_{a}\varphi_{\varepsilon}(y,z)\right|=0$
for all $0<\rho<N/4${\rm .}
\el

\vglue8pt

\demo{Proof}
(a) Throughout the proof $\varepsilon,N$ are fixed.
\medbreak
\hglue18pt {i.} We begin by removing both the bridge restriction and the torus
restriction.
For $y,z \in \Lambda_N$ and $0<b<N/2$, let
\be{qbdef}
q_b(y,z) = P_{y,z}(\sigma_b \leq \varepsilon),
\ee
where $\sigma$ is the first entrance time into $B_b$. There exists a constant
$c_2$ such that $\frac{dP_{y,z}}{dP_y} \leq c_2$, where $P_y(\cdot) =
P((\beta(t))_{t \in [0,\varepsilon/2]} \in \cdot \mid \beta(0)=y)$. Hence
\be{rembrd}
\sup\limits_{y,z \notin B_{b'}} q_b(y,z) \leq 2c_2
\sup\limits_{y \notin B_{b'}} P_y (\sigma_b \leq \varepsilon/2)
\qquad \mbox{ for all } 0<b<b'<N/4.\enspace
\ee
Let $\widehat\sigma_{N/2}$ be the first entrance time into $B^c_{N/2}
= \Lambda_N \setminus B_{N/2}$. Then for any $y \notin B_{b'}$ we may
decompose
\be{decom}
P_y(\sigma_b \leq \varepsilon/2) = P_y(\sigma_b \leq \varepsilon/2,\sigma_b
< \widehat \sigma_{N/2}) + P_y(\sigma_b \leq \varepsilon/2, \sigma_b
\geq \widehat\sigma_{N/2}). \quad
\ee
To estimate the second term on the right-hand side, we note that on its way
from $\partial B_{N/2}$ to $\partial B_b$ the Brownian motion must first
cross $\partial B_{N/4}$ and then cross $\partial B_{b'}$. Hence for any
$y \notin B_{b'}$,
\be{decom2}
P_y(\sigma_b \leq \varepsilon/2, \sigma_b \geq \widehat\sigma_{N/2})
\leq c_3 \sup\limits_{x \in \partial B_{b'}} P_x(\sigma_b \leq \varepsilon/2)
\ee
with $c_3 = \sup_{x \in \partial B_{N/2}} P_x(\sigma_{N/4} \leq \varepsilon/2)$.
Evidently, $c_3<1$ and so we deduce from (\ref{decom}) and (\ref{decom2}) that
\be{decom3}
\sup\limits_{y \notin B_{b'}} P_y (\sigma_b \leq \varepsilon/2)
\leq \frac{1}{1-c_3} \sup\limits_{y \notin B_{b'}}
P_y (\sigma_b \leq \varepsilon/2, \sigma_b < \widehat\sigma_{N/2}).\hskip.5in
\ee
As long as the Brownian motion does not hit $B^c_{N/2}$ it behaves like a
Brownian motion on $\R^d$. Therefore
\be{decom4}
P_y (\sigma_b \leq \varepsilon/2, \sigma_b < \widehat\sigma_{N/2})
\leq P_y^\infty(\sigma_b \leq \varepsilon/2),
\ee
where the upper index $\infty$ refers to removal of the torus restriction.
Combining (\ref{rembrd}), (\ref{decom3})  and (\ref{decom4}) we arrive at
\vglue-9pt
\be{decom5}
\sup\limits_{y,z \notin B_{b'}} q_b(y,z) \leq \frac{2c_2}{1-c_3}
\sup\limits_{y \notin B_{b'}} P_y^\infty(\sigma_b \leq \varepsilon/2)
\qquad \mbox{ for all } 0<b<b'<N/4.
\ee
\medbreak
 
\hglue18pt {ii.} Since
\be{hitest}
P_y^\infty(\sigma_b \leq \varepsilon/2)
\leq P_y^\infty(\sigma_b < \infty) = \left(\frac{b}{|y|}\right)^{d-2},
\ee
we obtain from (\ref{decom5}) that
\be{hitest2}
\sup\limits_{y,z \notin B_{b'}} q_b(y,z) \leq \frac{2c_2}{1-c_3}
\left(\frac{b}{b'}\right)^{d-2}.
\ee
Now put $b=a/\sqrt{T_\tau}$, $b'=K/\sqrt{T_\tau}$ and take the limit
$K \ra \infty$, to get the claim in Part (a).\medbreak

(b) Throughout the proof $\varepsilon,N$ are again fixed.
\smallbreak
\hglue18pt  i.  We shall prove that
\be{L2b1}
\lim_{b \da 0} \sup_{y,z \notin B_\rho} \left|\frac{q_b(y,z)}{\kappa_b}
- \varphi_\varepsilon(y,z)\right| = 0 \qquad \mbox{ for all } \rho>0.
\ee
Put $b=a/\sqrt{T_\tau}$ in (\ref{L2b1}) (recall (\ref{14})) and use the fact that
$\kappa_{a/\sqrt{T_\tau}}
= \kappa_a/\tau$ (recall (\ref{Tdef})) to get the claim in Part (b).

\medbreak \hglue18pt {ii.} Let $0<\delta<\varepsilon/2$. Define
\be{L2b2}
q_b^\delta(y,z) = P_{y,z}(\sigma_b \in [\delta,\varepsilon-\delta]).
\ee
Then, by the argument in Step i of Part (a),
\be{L2b3}
\ba{lll}
\sup\limits_{y,z \notin B_\rho} |q_b(y,z) - q_b^\delta(y,z)|
&\leq& \sup\limits_{y,z\notin B_\rho} P_{y,z}(\exists s \in [0,\delta]
\cap [\varepsilon-\delta,\varepsilon] \colon ~\beta(s) \in B_b) \\ \\
&\leq& \frac{2c_2}{1-c_3} \sup\limits_{y \notin B_\rho}
P_y^\infty(\sigma_b \leq \delta).
\ea
\ee
The supremum on the right-hand side is taken at any $y_0 \in \partial B_\rho$.
We may now invoke a result by Le Gall \cite{LG86}, which says that
\be{LG1}
\lim_{b \da 0} \frac{1}{\kappa_b} P_y^\infty(\sigma_b \leq t)
= \int_0^t p_s(-y) ds \qquad \mbox{ for all } y \in \R^d, t \geq 0. \hskip.5in
\ee
This gives us
\be{L2b4}
\lim_{b \da 0} \frac{1}{\kappa_b}
\sup\limits_{y,z \notin B_\rho} |q_b(y,z) - q_b^\delta(y,z)|
\leq  \frac{2c_2}{1-c_3} \int_0^{\delta} p_s(-y_0) ds. \hskip.5in
\ee
We thus see that to prove (\ref{L2b1}) it suffices to show that
\be{L2b5}
\lim_{\delta \da 0} \lim_{b \da 0} \sup_{y,z \notin B_\rho}
\left|\frac{q_b^\delta(y,z)}{\kappa_b}
- \varphi_\varepsilon(y,z)\right| = 0 \qquad \mbox{ for all } \rho>0. \hskip.5in
\ee

\bigbreak \hglue18pt {iii.} To analyze $q_b^\delta(y,z)$, we make a first entrance decomposition
on $\partial B_b$:
\be{L2b6}
q_b^\delta(y,z) = \frac{1}{p_{\varepsilon}(z-y)}
\int_\delta^{\varepsilon-\delta} \int_{\partial B_b} P_y(\sigma_b \in ds,
\beta(\sigma_b) \in dx) ~p_{\varepsilon-s}(z-x).\hskip.5in
\ee
Next we note that $p_{\varepsilon-s}(z-x)=[1+o_\delta(1)]p_{\varepsilon-s}(z)$
uniformly in $z \notin B_\rho$, $x \in \partial B_b$, $s \in
[\delta,\varepsilon-\delta]$,
where the $o_\delta(1)$ refers to $b \da 0$ for fixed $\delta$. Inserting this
approximation into (\ref{L2b6}), we get
\be{L2b7}
q_b^\delta(y,z) = \frac{1+o_\delta(1)}{p_{\varepsilon}(z-y)}
\int_\delta^{\varepsilon-\delta} P_y(\sigma_b \in ds) ~p_{\varepsilon-s}(z).
\ee
For the full integral we have
\be{L2b8}
\int_0^\varepsilon P_y(\sigma_b \in ds) ~p_{\varepsilon-s}(z)
= \int_0^\varepsilon ds ~P_y(\sigma_b \leq s) \left[\frac{\partial}{\partial s}
p_{\varepsilon-s}(z)\right] \hskip.5in
\ee
and so using (\ref{LG1}) we get, by dominated convergence,
\bel{L2b9}
\qquad\quad\lim\limits_{b \da 0} \frac{1}{\kappa_b}
\int_0^\varepsilon P_y(\sigma_b \in ds) ~p_{\varepsilon-s}(z)
&\hskip-6pt =\hskip-6pt&\int_0^\varepsilon ds \int_0^s du ~p_u(-y)
\left[\frac{\partial}{\partial s} p_{\varepsilon-s}(z)\right]\\
&\hskip-6pt =\hskip-6pt&\int_0^\varepsilon ds ~p_s(-y) p_{\varepsilon-s}(z).
\nne
The limit is in fact uniform in $y,z \notin B_\rho$, because $\Lambda_N$
is a compact set. Therefore, recalling (\ref{L2b7}) and the definition of
$\varphi_\varepsilon(y,z)$
in (\ref{15}), we see that to prove (\ref{L2b5}) it suffices to show that
uniformly in
$y,z \notin B_\rho$,
\be{L2b10}
\lim_{\delta \da 0} \lim_{b \da 0} \frac{1}{\kappa_b}
\int_0^\delta P_y(\sigma_b \in ds) ~p_{\varepsilon-s}(z) = 0,
\ee
and similarly for the integral over $[\varepsilon-\delta,\varepsilon]$.
However, the second factor is bounded uniformly in $z \notin B_\rho$ and $s
\in [0,\delta]$,
and so we are left with\break $\frac{1}{\kappa_b}P_y(\sigma_b \leq \delta)$.
Since $P_y(\sigma_b \leq \delta)
\leq \frac{1}{1-c_3}P^\infty_y(\sigma_b \leq \delta)$ for all $0<\rho<N/4$,
by the argument in Step
i of Part (a), we indeed get (\ref{L2b10}) via another application of
(\ref{LG1}).
\enddemo

  4.  We pick up the line of proof left off at the end of Step 2.
{}From Lemma \ref{cut1}(a) it follows that there exists $\delta_{K}>0$,
satisfying $\lim_{K\rightarrow\infty}\delta_{K}=0$, such that
\be{20}
-(1+\delta_{K})q_{\tau,\varepsilon}^{K/\sqrt{T_{\tau}}}
\leq \log\left(1-q_{\tau,\varepsilon}^{K/\sqrt{T_{\tau}}}\right)
\leq -q_{\tau,\varepsilon}^{K/\sqrt{T_{\tau}}}.
\ee
We are therefore naturally led to an investigation of the functions
$$\Phi_{\tau,\eta,\rho} \colon {\cal M}_{1}^{+}(\Lambda_{N}\times\Lambda_{N})
\mapsto [0,\infty)$$ defined by \vglue-9pt
\be{21}
\Phi_{\tau,\eta,\rho}(\mu) = \int_{\Lambda_{N}}dx
\left(1-\exp\left[-\eta\tau\int_{\Lambda_N\times\Lambda_N}
q_{\tau,\varepsilon}^{\rho}(y-x,z-x)
\mu(dy,dz)\right]\right),
\ee
for which (\ref{repr}) and (\ref{20}) give us the following sandwich:
\be{sand*}
\Phi_{\tau,1/\varepsilon,K/\sqrt{T_{\tau}}}(L_{\tau,\varepsilon})
\leq \E_{\tau,\varepsilon}(V_{\tau}^{K})
\leq \Phi_{\tau,(1+\delta_{K})/\varepsilon,K/\sqrt{T_{\tau}}} (L_{\tau,\varepsilon}). \hskip.5in
\ee
The functions $\Phi_{\tau,\eta,\rho}$ have nice continuity properties:

\specialnumber{3} \proclaim{Lemma} \label{cont1}
There exist constants $c_4,c_5$ such that\/{\rm :}
\begin{itemize}
\item[{\rm (a)}] $\left|\Phi_{\tau,\eta,\rho}(\mu)-\Phi_{\tau,\eta,\rho^{\prime}}(\mu)\right|
\leq c_4\eta\left|\sqrt{\rho}-\sqrt{\rho^{\prime}}\right|$ for all
$\eta,\mu$ and
$\tau\geq\tau_0(\rho,\rho^{\prime})${\rm .}

\item[{\rm (b)}] $\left|\Phi_{\tau,\eta,\rho}(\mu)-\Phi_{\tau,\eta^{\prime},\rho}(\mu)\right|
\leq c_5\left|\eta-\eta^{\prime}\right|$ for all $\rho,\mu$ and $\tau\geq\tau_0(\rho)${\rm .}
\end{itemize}

\endproclaim

\demo{Proof}
The proof uses Lemma \ref{cut1}(b). 
\smallbreak
(a) Write
\begin{eqnarray}\label{lcontaproof}
&&\\
&&\hskip-12pt\left|\Phi_{\tau,\eta,\rho}(\mu)-\Phi_{\tau,\eta,\rho^{\prime}}(\mu)\right|
\nn\\
&&\qquad \leq \eta \int_{\Lambda_N} dx \int_{\Lambda_N\times\Lambda_N} \mu(dy,dz)
\left|\tau q^{\rho}_{\tau,\varepsilon}(y-x,z-x)
-\tau q^{\rho^{\prime}}_{\tau,\varepsilon}(y-x,z-x)\right|\nonumber \\
&&\qquad = \eta \int_{\Lambda_N} dx \int_{\Lambda_N\times\Lambda_N} \mu(dy,dz)
\left[\left|\kappa_a\varphi_{\varepsilon}^{\rho}(y-x,z-x)\phantom{\varphi_{\varepsilon}^{\rho^{\prime}}}\right.\right.\nn\\
&&\hskip2.35in\left.\left.
-\kappa_a\varphi_{\varepsilon}^{\rho^{\prime}}(y-x,z-x)\right|
+ o_{\rho,\rho^{\prime}}(1)\right]\nonumber \\
&&\qquad = \eta \left[\left|\sqrt{\rho}-\sqrt{\rho^{\prime}}\right|
+ \left|\Lambda_N\right|o_{\rho,\rho^{\prime}}(1)\right].
\nonumber
\end{eqnarray}
Here, $o_{\rho,\rho^{\prime}}(1)$ means an error tending to zero as
$\tau\ra\infty$
depending on $\rho,\rho^{\prime}$, and in the last equality we use the fact that
$\int_{\Lambda_N}dx\,
\varphi_{\varepsilon}(y-x,z-x)=\varepsilon$ for\break all $y,z$. 
\medbreak
(b) Write
\begin{eqnarray} \label{lcont1bproof}
&&\\
 &&\hskip-.5in\left|\Phi_{\tau,\eta,\rho}(\mu)-\Phi_{\tau,\eta^{\prime},\rho}(\mu)\right|
 \nn \\[4pt]
&&\qquad \leq \left|\eta-\eta^{\prime}\right| \int_{\Lambda_N} dx
\int_{\Lambda_N\times\Lambda_N} \mu(dy,dz)\,
\tau q^{\rho}_{\tau,\varepsilon}(y-x,z-x)\nonumber \\[4pt]
&&\qquad = \left|\eta-\eta^{\prime}\right| \int_{\Lambda_N} dx
\int_{\Lambda_N\times\Lambda_N} \mu(dy,dz)
\left[\kappa_a\varphi_{\varepsilon}^{\rho}(y-x,z-x)+o_{\rho}(1)\right]\nonumber \\[4pt]
&&\qquad \leq  \left|\eta-\eta^{\prime}\right| \left[\kappa_a \varepsilon +
\left|\Lambda_N\right|o_{\rho}(1)\right],\nonumber
\end{eqnarray}
where in the last inequality we drop the superscript $\rho$ to be able to
perform the $x$-integration.
\enddemo

 5.  With the help of (\ref{cut}), (\ref{sand*}) and Lemma
\ref{cont1}(a), (b),
we get
\begin{eqnarray} \label{22}
&&\\
\E_{\tau,\varepsilon}(V_{\tau})
&\leq& \E_{\tau,\varepsilon}(V_{\tau}^{K})+c_1K^{d}/T_{\tau}\nonumber \\[4pt]
&\leq& \Phi_{\tau,(1+\delta_{K})/\varepsilon,K/\sqrt{T_{\tau}}}(L_{\tau,\varepsilon})
+c_1K^{d}/T_{\tau}\nonumber \\[4pt]
&\leq& \Phi_{\tau,1/\varepsilon,\rho}(L_{\tau,\varepsilon})+ c_1K^{d}/T_{\tau}
+ c_4\left[\sqrt{K/\sqrt{T_{\tau}}}+\sqrt{\rho}\right]/\varepsilon + c_5\delta_{K}/\varepsilon,
\nonumber
\end{eqnarray}
and also a similar lower bound.
\medbreak
 6.  Next we approximate
$\Phi_{\tau,1/\varepsilon,\rho}(L_{\tau,\varepsilon})$ by
$\Phi_{\infty,1/\varepsilon,\rho}(L_{\tau,\varepsilon})$ defined as \vglue-9pt
\be{23}
\Phi_{\infty,\eta,\rho}(\mu)=\int_{\Lambda_{N}}dx
\left(1-\exp\left[-\eta\kappa_{a}\int_{\Lambda_N\times\Lambda_N}
\varphi_{\varepsilon}^{\rho}(y-x,z-x)\mu(dy,dz)\right]\right),
\ee
where for $\rho>0$ we define
$\varphi_{\varepsilon}^{\rho}(y,z)=\varphi_{\varepsilon}(y,z)$ if
$y,z\notin B_{\rho}$ and zero otherwise. For that we need the following:

\specialnumber{4} \proclaim{Lemma} \label{cont}
There exist constants $c_6,c_7>0$ such that\/{\rm :}
\begin{itemize}
\item[{\rm (a)}] $\left|\Phi_{\infty,\eta,\rho}(\mu)-\Phi_{\tau,\eta,\rho}(\mu)\right|\leq
c_6\eta\delta_{\rho,\tau}$ for all $\mu$ with
$\lim_{\tau\ra\infty}\delta_{\rho,\tau}
=0$ for any $\rho>0${\rm .}

\item[{\rm (b)}] $\left|\Phi_{\infty,1/\varepsilon,0}(\mu)-\Phi_{\infty,1/\varepsilon,0}
(\mu^{\prime})\right|\leq c_7\left\|\mu-\mu^{\prime}\right\|_{tv}${\rm ,} where
$\|\cdot\|_{tv}$ denotes the total variation norm{\rm .}$\phantom{\sum^\int}$\end{itemize}

\el

\demo{Proof}
The proof again uses Lemma \ref{cut1}(b).

(a) Write
\begin{eqnarray}\label{lconta}
&&\\
&&\hskip-12pt \left|\Phi_{\infty,\eta,\rho}(\mu)-\Phi_{\tau,\eta,\rho}(\mu)\right| \nonumber \\[4pt]
&&\qquad \leq \eta \int_{\Lambda_N} dx \int_{\Lambda_N\times\Lambda_N} \mu(dy,dz)
\left[\tau q_{\tau,\varepsilon}^{\rho}(y-x,z-x)
-\kappa_a\varphi_{\varepsilon}(y-x,z-x)\right]\nonumber \\[4pt]
&&\qquad = \eta \int_{\Lambda_N} dx \int_{\Lambda_N\times\Lambda_N} \mu(dy,dz)
\,o_{\rho}(1) = \eta \left|\Lambda_N\right| o_{\rho}(1).
\nonumber
\end{eqnarray}

(b) Write
\begin{eqnarray}\label{lcontb}
&&\left|\Phi_{\infty,1/\varepsilon,0}(\mu)-
\Phi_{\infty,1/\varepsilon,0}(\mu^{\prime})
\right|  \\[4pt]
&&\qquad \leq \frac{\kappa_a}{\varepsilon} \int_{\Lambda_N} dx
\int_{\Lambda_N\times\Lambda_N}
\left|\mu-\mu^{\prime}\right|(dy,dz) \,\varphi_{\varepsilon}(y-x,z-x)\nonumber \\[4pt]
&&\qquad = \kappa_a \int_{\Lambda_N\times\Lambda_N}
\left|\mu-\mu^{\prime}\right|(dy,dz)
= \kappa_a  \left\|\mu-\mu^{\prime}\right\|_{tv}.\nonumber
\end{eqnarray}
\vglue-24pt
\enddemo

\phantom{snow}

  7.  Using (\ref{22}), the similar lower bound and Lemma \ref{cont}(a) with
$\eta=1/\varepsilon$, we now have that for any $K$ and $\rho$,
\bel{25}
&&\left\|\E_{\tau,\varepsilon}(V_{\tau})-\Phi_{\infty,1/\varepsilon,0}
(L_{\tau,\varepsilon})\right\|_{\infty}\\[4pt]
&&\qquad\qquad
\leq c_1 K^{d}/T_{\tau} + c_4 \left[\sqrt{K/\sqrt{T_{\tau}}}
+\sqrt{\rho}\right]/\varepsilon + c_5 \delta_{K}/\varepsilon
+c_6\delta_{\rho,\tau}/\varepsilon.
\nne
Letting $\tau\ra\infty$, followed by $K \ra \infty$ and $\rho \da 0$, we thus
arrive at 
\be{26}
\lim_{\tau\rightarrow\infty}\left\|
\E_{\tau,\varepsilon}(V_{\tau})-\Phi_{\infty,1/\varepsilon,0}
(L_{\tau,\varepsilon})\right\|_{\infty}=0
\qquad \mbox{ for all } \varepsilon>0.\hskip.5in
\ee

 8.  The desired LDP for fixed $\varepsilon$ can now be derived as follows.
First, note that $\Phi_{\infty,1/\varepsilon,0}$ is continuous by Lemma
\ref{cont}(b)
(even in the total variation topology). Next, note from (\ref{23}) that
$\Phi_{\infty,1/\varepsilon,0}= \Phi_{1/\varepsilon}$, the function 
defined in
(\ref{Phiext}). Therefore we can use one of the standard results of Donsker and
Varadhan \cite{DV*}(III) (see also Bolthausen \cite{Bo2}), namely, that
$(L_{\tau,\varepsilon})_{\tau>0}$ satisfies the LDP on
${\cal M}^+_1(\Lambda_N\times\Lambda_N)$ with rate $\tau$ and with rate
function $\frac{1}{\varepsilon}I_{\varepsilon}^{(2)}$ defined in (\ref{27}).
{}From the contraction principle and (\ref{26}) we now get the claim
in Proposition \ref{Epsilon_LDP}.\hfill\qed
\enddemo

\demo{{\rm 2.4.} The limit $\varepsilon\da 0$ for the {\rm LDP}}
We already know from Section 2.2 that the quantity of interest, namely
$V_{\tau}$, is for small $\varepsilon$ well approximated by
$\E_{\tau,\varepsilon}(V_{\tau})$,
for which we have the LDP in Proposition \ref{Epsilon_LDP}. The main step
to prove the LDP
for $V_{\tau}$ itself is therefore to derive an appropriate limit result
for the rate function
in (\ref{28}). This needs some preparations.

\medbreak
{1.} We denote by $I \colon {\cal M}_{1}^{+}(\Lambda_{N})\mapsto
[0,\infty]$
the standard large deviation rate function for the empirical distribution of
the Brownian motion:
\be{29}
\begin{array}{llll}
I(\nu)&=&\frac{1}{2}\int_{\Lambda_N}|\nabla\phi|^2(x)dx
&\mbox{if $\frac{d\nu}{dx}=\phi^2$ with $\phi \in H^1(\Lambda_N)$}\\[4pt]
&=&\infty &\mbox{otherwise}.
\end{array}\hskip.5in
\ee
We further denote by $I_{\varepsilon} \colon {\cal M}_{1}^{+}(\Lambda_{N})
\mapsto [0,\infty]$ the following projection of $I_{\varepsilon}^{(2)}$
(recall (\ref{27}))
onto ${\cal M}_{1}^{+}(\Lambda_{N})$:
\be{30}
I_{\varepsilon}(\nu)=\inf\left\{I_{\varepsilon}^{(2)}(\mu) \colon
\,\mu_{1}=\nu\right\}.
\ee
We begin by collecting some basic facts about these entropies, all of which
have
been proved by Donsker and Varadhan \cite{DV*} or are simple consequences
of their results:

\specialnumber{5} \proclaim{Lemma} \label{Entropy}
Let $(\pi_t)_{t\geq 0}$ denote the semigroup of the Brownian motion{\rm .} Then
for all $\nu,\mu$\/{\rm :}\/
\begin{itemize}
\ritem{(a)} $I_t(\nu)=-\inf_{u\in{\cal D}^{+}}
\int\log\frac{\pi_{t}u}{u}d\nu$, where ${\cal D}^{+}$ is
the set of positive measurable functions bounded away from $0$ and $\infty${\rm .}

\ritem{(b)} $t\mapsto I_t(\nu)/t$ is nonincreasing with $\lim_{t\da 0}
I_t(\nu)/t = I(\nu)${\rm .}

\ritem{(c)} $\left\|\nu-\nu \pi_{s}\right\|_{tv} \leq 8\sqrt{I_{s}(\nu)}$ for $s>0${\rm .} 
\ritem{(d)} $I_s(\nu \pi_{t})\leq I_s(\nu)$ for $s,t>0${\rm .}
\ritem{(e)} $\left\|\mu-\mu_{1}\otimes \pi_{s}\right\|_{tv} \leq
8\sqrt{I_{s}^{(2)}(\mu)}$
for $s>0${\rm .}
\end{itemize}

\el

\demo{Proof}
(a) This is \cite[(III), Th.~2.1]{DV*}, combined with \cite[(I),
Lemma 2.1]{DV*}.
\medbreak
(b) Fix $s,t>0$. For every $u\in{\cal D}^{+}$,
\be{31}
\int\log\frac{\pi_{s+t}u}{u}d\nu=\int\log\frac{\pi_{s}(\pi_{t}u)}{\pi_{t}u}d\nu
+\int\log\frac{\pi_{t}u}{u}d\nu\geq-I_{s}(\nu)-I_{t}(\nu).\quad
\ee
Taking the infimum over $u$ and using (a), we get $-I_{s+t}(\nu)\geq
-I_{s}(\nu)-I_{t}(\nu)$.
Hence $t \mapsto I_{t}(\nu)/t$ is nonincreasing. The fact that $\lim_{t\da 0}
I_{t}(\nu)/t=I(\nu)$ is \cite[(I), Lemma 3.1]{DV*}.
\medbreak
(c) This is \cite[(I), Lemma 4.1]{DV*}.
\medbreak
(d) This follows from the convexity of $\nu \mapsto I_{s}(\nu)$ for $s>0$.
\medbreak
(e) Let $P^{\mu}(x,dy)$ be any transition kernel on $\Lambda_{N}$ such that
$\mu=\mu_{1}\otimes P^{\mu}${\rm .} Then
\be{32}
\left\|\mu-\mu_{1}\otimes \pi_{s}\right\|_{tv}
\leq\int\mu_{1}(dx)\left\|
P^{\mu}(x,\cdot)-\pi_{s}(x,\cdot)\right\|_{tv}.
\ee
By Csisz\'ar \cite[Th.\ 4.1]{C}, we have (recall that $h(\cdot|\cdot)$
denotes
relative entropy)
\be{33}
\left\|P^{\mu}(x,\cdot)-\pi_{s}(x,\cdot)\right\|_{tv}
\leq 8 \sqrt{h(P^{\mu}(x,\cdot)|\pi_{s}(x,\cdot))}.
\ee
Therefore
\bel{34}
\qquad\quad\left\|\mu-\mu_{1}\otimes \pi_{s}\right\|_{tv}
&\leq& 8\int\mu_{1}(dx)\sqrt{h(P^{\mu}(x,\cdot)|\pi_{s}(x,\cdot))}  \\[4pt]
&\leq& 8\sqrt{\int\mu_{1}(dx)h(P^{\mu}(x,\cdot)|\pi_{s}(x,\cdot))}
=8\sqrt{I_{s}^{(2)}(\mu)},
\nne
where the last equality uses (\ref{27}).
\enddemo

  2.  To take advantage of the link provided by Lemma \ref{Entropy}(b),
we shall need an approximation of the functions $\Phi_{1/\varepsilon}
\colon {\cal M}_{1}^{+}(\Lambda_{N}\times \Lambda_{N})\mapsto [0,\infty)$,
appearing in Proposition \ref{Epsilon_LDP}, by the simpler functions
$\Psi_{1/\varepsilon} \colon {\cal M}_{1}^{+}(\Lambda_{N})\mapsto [0,\infty)$
defined by
\be{35}
\Psi_{1/\varepsilon}(\nu)=\int_{\Lambda_N} dx
\left(1-\exp\left[-\frac{\kappa_{a}}{\varepsilon}
\int_{0}^{\varepsilon} ds \int_{\Lambda_N} p_{s}(x-y)\nu(dy)\right]\right).\hskip.25in
\ee

\phantom{snow}
\vglue-6pt
\specialnumber{6} \proclaim{Lemma} \label{Approx1}
For any $K>0${\rm ,}
\be{36}
\lim_{\varepsilon\da 0}\sup_{\mu\colon~\frac{1}{\varepsilon}
I_{\varepsilon}^{(2)}(\mu)\leq K}\left|\Phi_{1/\varepsilon}
(\mu)-\Psi_{1/\varepsilon}(\mu_{1})\right| = 0.
\ee
\el

\demo{Proof}
As is obvious by comparing (\ref{Phiext}), (\ref{15}) with (\ref{35}),
we have $\Psi_{1/\varepsilon}(\mu_{1})\break=\Phi_{1/\varepsilon}
(\mu_{1}\otimes \pi_{\varepsilon}).$  Therefore
\begin{eqnarray} \label{37}
&&\\
&&\hskip-12pt \left|\Phi_{1/\varepsilon}(\mu)-\Psi_{1/\varepsilon}(\mu_{1})\right|\nn  \\[4pt]
&&\qquad = \left|\Phi_{1/\varepsilon}(\mu)-\Phi_{1/\varepsilon}
(\mu_{1}\otimes \pi_{\varepsilon})\right|\nn \\[4pt]
&&\qquad \leq \frac{\kappa_{a}}{\varepsilon}\left|\int_{\Lambda_N}
dx\int_{\Lambda_{N}
\times\Lambda_{N}}\varphi_{\varepsilon}(y-x,z-x)\left[\mu(dy,dz)-(\mu_{1}
\otimes \pi_{\varepsilon})(dy,dz)\right]\right|\nn \\[4pt]
&&\qquad \leq \frac{\kappa_{a}}{\varepsilon} \int_{\Lambda_{N}\times\Lambda_{N}}
\left\{\int_{\Lambda_N} dx\,\varphi_{\varepsilon}(y-x,z-x)\right\}
\left|\mu-\mu_{1}\otimes \pi_{\varepsilon}\right|(dy,dz)\nn \\[4pt]
&&\qquad = \kappa_{a}\left\|\mu-\mu_{1}\otimes \pi_{\varepsilon}\right\|_{tv},
\nn
\end{eqnarray}
where in the last equality we again use the fact  that the integral between braces
equals $\varepsilon$
for all $y,z$ by (\ref{15}). The claim now follows from Lemma \ref{Entropy}(e).
\enddemo

  3.  Next, we define the function $\Gamma\colon
L_{1}^{+}(\Lambda_{N})\mapsto
[0,\infty)$ by
\be{38}
\Gamma(f)=\int_{\Lambda_N} dx \left(1-e^{-\kappa_{a}f(x)}\right).
\ee

\specialnumber{7} \proclaim{Lemma} \label{Approx2}
For any $K>0,$
\be{39}
\lim_{\varepsilon\da 0}\sup_{\nu \colon ~\frac{1}{\varepsilon}
I_{\varepsilon}(\nu)\leq K}\left|\Gamma\left(\textstyle{\frac{d\nu}{dx}}\right)
-\Psi_{1/\varepsilon}(\nu)\right| = 0.
\ee
{\rm (}\/Note from {\rm (\ref{27})} and {\rm (\ref{30})} that if $I_{\varepsilon}(\nu)<\infty${\rm ,} then
$d\nu\ll dx$ because $\nu\otimes\pi_{\varepsilon}\ll dx\otimes dy${\rm .)}
\el

\demo{Proof}
Using (\ref{35}) and (\ref{38}), we write
\begin{eqnarray}\label{40}
&&\left|\Gamma\left(\frac{d\nu}{dx}\right)-\Psi_{1/\varepsilon}(\nu)\right|  \\[4pt]
&&\qquad \leq \int_{\Lambda_N} dx\left|\exp\left[-\frac{\kappa_{a}}{\varepsilon}
\int_{0}^{\varepsilon}ds \int_{\Lambda_N} p_{s}(x-y)\nu(dy)\right]\right.
\nn\\[4pt]
&&\hskip1in\left.
-\exp\left[-\frac{\kappa_{a}}{\varepsilon}\int_{0}^{\varepsilon}ds\,
\frac{d\nu}{dx}(x)\right]\right|\nn \\[4pt]
&&\qquad \leq\int_{\Lambda_N}
dx\,\frac{\kappa_{a}}{\varepsilon}\int_{0}^{\varepsilon}ds
\left|\frac{\nu \pi_{s}}{dx}(x)-\frac{d\nu}{dx}(x)\right|
=\frac{\kappa_{a}}{\varepsilon}\int_{0}^{\varepsilon}ds\left\|\nu
\pi_{s}-\nu\right\|_{tv}.
\nn
\end{eqnarray}
Now, for $0\leq s\leq\varepsilon$ we have, by Lemma \ref{Entropy}(c),
\be{41}
\left\|\nu \pi_{s}-\nu\right\|_{tv} \leq
\left\|\nu \pi_{s}-\nu \pi_{s+\varepsilon}\right\|_{tv}
+\left\|\nu \pi_{s+\varepsilon}-\nu\right\|_{tv}
\leq 8 \sqrt{I_{\varepsilon}(\nu \pi_{s})}+8\sqrt{I_{\varepsilon+s}(\nu)}.
\ee
Moreover, $I_{\varepsilon}(\nu \pi_{s})\leq I_{\varepsilon}(\nu)$ by Lemma
\ref{Entropy}(d) and $$I_{\varepsilon+s}(\nu)\leq 2\varepsilon I_{\varepsilon+s}
(\nu)/(\varepsilon+s)\leq 2\varepsilon I_{\varepsilon}(\nu)/\varepsilon = 2
I_{\varepsilon}(\nu)$$
by Lemma \ref{Entropy}(b). Thus we get $\left\|\nu \pi_{s}-\nu\right\|_{tv}
\leq 8(1+\sqrt{2})\sqrt{I_{\varepsilon}(\nu)}$. From this the claim follows.
\enddemo

We now have all the ingredients to perform the proof of Proposition \ref{VLDP}
in Section 2.1.

\demo{{\rm 2.5.} Proof of Proposition {\rm \ref{VLDP}}}
For any $f \colon \R_+ \mapsto \R$ bounded and continuous:
\begin{eqnarray}\label{42}
&&\\
\lim\limits_{\tau\rightarrow\infty}\frac{1}{\tau}\log
E\left(\exp\left[\tau f(V_{\tau})\right]\right)
&\hskip-8pt =\hskip-6pt & \lim\limits_{\varepsilon\da 0}\lim\limits_{\tau\rightarrow\infty}
\frac{1}{\tau}\log E\left(\exp\left[\tau f(\E_{\tau,\varepsilon}
(V_{\tau}))\right]\right) \nn \\
&\hskip-8pt =\hskip-6pt &\lim\limits_{\varepsilon\da 0}\sup\limits_{\mu}\left
\{f(\Phi_{1/\varepsilon}(\mu))-\frac{1}{\varepsilon}
I_{\varepsilon}^{(2)}(\mu)\right\}\nn \\
&\hskip-8pt =\hskip-6pt &\lim\limits_{K\rightarrow\infty}\lim\limits_{\varepsilon\da 0}
 \sup\limits_{\mu \colon \,\frac{1}{\varepsilon}I_{\varepsilon}^{(2)}(\mu)\leq K}\left
\{f(\Phi_{1/\varepsilon}(\mu))-\frac{1}{\varepsilon}
I_{\varepsilon}^{(2)}(\mu)\right\}\nn \\
&\hskip-8pt =\hskip-6pt &\lim\limits_{K\rightarrow\infty}\lim\limits_{\varepsilon\da 0}
 \sup\limits_{\mu \colon  \frac{1}{\varepsilon}I_{\varepsilon}^{(2)}(\mu)\leq K}\left
\{f(\Psi_{1/\varepsilon}(\mu_{1}))-\frac{1}{\varepsilon}
I_{\varepsilon}^{(2)}(\mu)\right\}\nn \\
&\hskip-8pt =\hskip-6pt &\lim\limits_{K\rightarrow\infty}\lim\limits_{\varepsilon\da 0}
\sup\limits_{\nu \colon  \frac{1}{\varepsilon}I_{\varepsilon}(\nu)\leq K}\left
\{f(\Psi_{1/\varepsilon}(\nu))-\frac{1}{\varepsilon}I_{\varepsilon}(\nu)\right\}\nn \\
&\hskip-8pt =\hskip-6pt &\lim\limits_{K\rightarrow\infty}\lim\limits_{\varepsilon\da 0}
~\sup\limits_{\nu \colon  \frac{1}{\varepsilon}I_{\varepsilon}(\nu)\leq K}\left
\{f\left(\Gamma\left(\frac{d\nu}{dx}\right)\right)
-\frac{1}{\varepsilon}I_{\varepsilon}(\nu)\right\}\nn \\
&\hskip-8pt =\hskip-6pt &\sup\limits_{\nu}\left\{f\left(\Gamma\left(
\frac{d\nu}{dx}\right)\right)-I(\nu)\right\}\nn \\
&\hskip-8pt =\hskip-6pt &\sup\limits_{\phi \in H^1(\Lambda_N) \colon \,\|\phi\|^2_2=1}
\left\{f(\phi^2)-\frac{1}{2}\|\nabla\phi\|^2_2\right\}.
\nn
\end{eqnarray}
Here we use, respectively, Proposition \ref{Concentration}, Proposition
\ref{Epsilon_LDP},
Lemma \ref{Approx1}, equation (\ref{30}), Lemma \ref{Approx2}, Lemma
\ref{Entropy}(b) and
equation (\ref{29}). Recalling (\ref{38}), we see that the claim now
follows from
the inverse of Varadhan's lemma proved in Bryc \cite{Br}.
\hfill\qed\enddemo

\demo{{\rm 2.6.} Proof of Proposition {\rm \ref{p2}}}
Throughout this section, $a>0$ and $b>0$ are fixed. For ease of notation we
introduce the following abbreviations:
\be{eqaH}
A(\phi)=\displaystyle\int_{\R^d}\phi^2(x)\,dx,\quad
B(\phi)=\displaystyle\int_{\R^d}(1-e^{-\kappa_a\phi^2(x)})\,dx,\quad
C(\phi)=\displaystyle\int_{\R^d}|\nabla\phi|^2(x)\,dx
\ee
for $\phi\in H^1(\R^d)$, and their counterparts $A_N(\phi_N)$,
$B_N(\phi_N)$, $C_N(\phi_N)$ for
$\phi_N\in H^1(\Lambda_N)$ with $\Lambda_N=[-\frac{N}{2},\frac{N}{2})^d$,
the $N$--torus with periodic boundary conditions.

\medbreak 
1.  `$I_N^{\kappa_a}(b)\le I^{\kappa_a}(b)$ for all $N>0$'.\\

For $\phi\in H^1(\R^d)$, let $\sigma_N\phi\in H^1(\Lambda_N)$ be defined by
\be{eqbH}
(\sigma_N\phi)^2(x)=\left\{\ba{ll}\sum_{k\in\Z^d}\phi^2(x+kN)&\quad(x\in\Lambda_
N)\\
0&\quad(x\notin \Lambda_N).
\ea\right.
\ee
Then
\be{eqcH}
A_N(\sigma_N\phi)=A(\phi),\quad
B_N(\sigma_N\phi)\le B(\phi),\quad
C_N(\sigma_N\phi)\le C(\phi),\hskip.5in
\ee
where the second and third statements hold because
$1-e^{-(f+g)}\le(1-e^{-f})+(1-e^{-g})$,
respectively, $(\nabla\sqrt{f^2+g^2})^2\le(\nabla f)^2+(\nabla g)^2$ for
arbitrary functions
$f,g\ge0$. Hence
\begin{eqnarray} \label{eqdH}
&&\\
I_N^{\kappa_a}(b)&=&\inf\{\,C_N(\phi_N) \colon ~\phi_N\in H^1(\Lambda_N),
A_N(\phi_N)=1,B_N(\phi_N)\le b\,\}\nn  \\
&=&\inf\{\,C_N(\sigma_N\phi) \colon ~\phi\in H^1(\R^d),
A_N(\sigma_N\phi)=1,B_N(\sigma_N\phi)\le b\,\}\nn \\
&\le&\inf\{\,C (\phi) \colon ~\phi\in H^1(\R^d),
A(\phi)=1,B(\phi)\le b\,\} = I^{\kappa_a}(b).
\nn
\end{eqnarray}

 2.  `$\liminf_{N\to\infty}I_N^{\kappa_a}(b)\ge I^{\kappa_a}(b+)$'.\\

For every $\varepsilon>0$ there exists a $\phi_N\in H^1(\Lambda_N)$ such that
\be{eqeH}
A_N(\phi_N)=1,\quad
B_N(\phi_N)\le b,\quad
C_N(\phi_N)\le I_N^{\kappa_a}(b)+\varepsilon,\hskip.5in
\ee
i.e., $\phi_N$ is an $\varepsilon$--minimiser. By shifting $\Lambda_N$
around, we see that there must exist a $y\in \Lambda_N$ such that
\be{eqfH}
\int_{\delta\Lambda_N}[\phi_N^2(x+y)+(\nabla\phi_N)^2(x+y)]\,dx
\le\frac{|\delta\Lambda_N|}{|\Lambda_N|}[A_N(\phi_N)+C_N(\phi_N)].
\ee
Let $\tau\phi_N\in H^1(\R^d)$ be defined by
\be{eqgH}
(\tau\phi_N)(x)=\left\{
\ba{ll}\phi_N(x+y)&\quad(x\in\Lambda_N)\\
\phi_N([x]_N+y)\Bigl\{N\Bigl(1-\frac{|x|}{|[x]_N|}\Bigr)+1\Bigr\}
&\quad(x\in\Lambda_{N+1}\setminus\Lambda_N)\\
0&\quad(x\notin\Lambda_{N+1})
\ea\right.
\ee
with $[x]_N$ the radial projection of $x$ onto $\delta\Lambda_N$, i.e.,
$\tau\phi_N$
linearly drops to 0 outside $\Lambda_N$ along radial lines. Then, clearly,
$(\tau\phi_N)^2(x)\le\phi_N^2([x]_N+y)$ and $(\nabla\tau\phi_N)^2(x)\le
d(\nabla\phi_N)^2([x]_N+y)$ for all $x\in \Lambda_{N+1}\setminus\Lambda_N$.
Hence, by (\ref{eqfH}),  
\bel{eqhH}
&&A(\tau\phi_N)\le A_N(\phi_N)+\delta_N,\quad
B(\tau\phi_N)\le B_N(\phi_N)+\kappa_a\delta_N,\\[4pt]
&&C(\tau\phi_N)  \le  C_N(\phi_N)+\delta_N
\nne
with
\be{eqiH}
\delta_N=d\frac{|\delta\Lambda_{N+1}|}{|\delta\Lambda_N|}\,\frac{|\delta\Lambda_
N|}{|\Lambda_N|}
[A_N(\phi_N)+C_N(\phi_N)]=O\biggl(\frac 1N\biggr).\hskip.5in
\ee
Now define $\phi^*\in H^1(\R^d)$ by
\be{eqjH}
\phi^*=\frac{\tau\phi_N}{\sqrt{A(\tau\phi_N)}}.
\ee
Then clearly
\be{eqkH}
A(\phi^*)=1,\quad
B(\phi^*)\le B(\tau\phi_N),\quad
C(\phi^*)=A(\tau\phi_N)C(\tau\phi_N),\hskip.5in
\ee
where the second statement holds because $A(\tau\phi_N)\geq A_N(\phi)=1$.
Combining (\ref{eqeH}), (\ref{eqhH}) and (\ref{eqkH}), we get
\be{eqlH}
A(\phi^*)=1,\quad
B(\phi^*)\le b+\kappa_a\delta_N,\quad
C(\phi^*)\le(1+\delta_N)
[I_N^{\kappa_a}(b)+\varepsilon+\delta_N].
\ee
Hence we have
\be{eqmH}
\ba{lll}
I^{\kappa_a}(b+\kappa_a\delta_N)&=&\inf\{\,C(\phi) \colon ~\phi\in
H^1(\R^d),A(\phi)=1,B(\phi)\le b+\kappa_a\delta_N\,\}\\ \\
&\le& C(\phi^*) \le (1+\delta_N)[I_N^{\kappa_a}(b)+\varepsilon+\delta_N].
\ea
\ee
Let $N\to\infty$ and use (\ref{eqiH}) to get
$I^{\kappa_a}(b+)\le\varepsilon+\liminf_{N\to\infty}I_N^{\kappa_a}(b)$. Since
$\varepsilon>0$
is arbitrary this proves the claim.\medbreak

  3.  Combining Steps 1 and 2 and noting that $b\mapsto I^{\kappa_a}(b)$
is right-continuous
(because it is nonincreasing and lower semicontinuous), we have completed
the proof
of Proposition \ref{p2}. \hfill\qed\enddemo

\section{The lower bound in Theorem \ref{t1}}

In this section we prove the complement of Propositions \ref{p1} and \ref{p2},\break
which will complete the proof of Theorem \ref{t1}. Recall from Section 2.1
that, by\break Brownian scaling, $t^{-1}|W^a(t)|$ has the same distribution as
$|W^{a\tau^{-1/(d-2)}}$ $(\tau)|$ with\break $\tau = t^{(d-2)/d}$.

\specialnumber{6} \proclaim{Proposition} \label{plb}
Let $d \geq 3$ and $a>0${\rm .} For every $b>0,$
\be{lb1}
\liminf\limits_{\tau \ra \infty} \frac{1}{\tau} \log
P(|W^{a\tau^{-1/(d-2)}}(\tau)| \leq b)
\geq - I^{\kappa_a}(b),
\ee
where $I^{\kappa_a}(b)$ is given by {\rm (\ref{mdsrf}),} {\rm (\ref{mdsset}).}
\ep

\demo{Proof}
Let $C_N(\tau)$ be the event that the Brownian motion does not hit
$\partial \Lambda_{N-a}$ until time $\tau$. Clearly,
\be{lb5}
P(|W^{a\tau^{-1/(d-2)}}(\tau)| \leq b)
\geq P\Big(C_N(\tau),\,|W_N^{a\tau^{-1/(d-2)}}(\tau)| \leq b\Big).\hskip.5in
\ee
The right-hand side involves the Brownian motion on the torus, but
restricted to stay a distance $a$ away from the boundary. We
can now simply repeat the argument in Section 2 on the event $C_N(\tau)$,
the result being that
\be{lb6}
\lim_{\tau\ra\infty} \frac{1}{\tau} \log
P\Big(|W_N^{a\tau^{-1/(d-2)}}(\tau)| \leq b ~|~ C_N(\tau)\Big)
= - \widetilde I_N^{\kappa_a}(b)
\ee
where $\widetilde I_N^{\kappa_a}(b)$ is given by the same formulas as in
(\ref{mdsrf}), (\ref{mdsset}), except that $\R^d$ is replaced by $\Lambda_N$
and $\phi$ is restricted to ${\rm supp}(\phi)\cap\partial\Lambda_N=\emptyset$.
We have
\be{addlambda}
\lim_{\tau\ra\infty} \frac{1}{\tau} \log P(C_N(\tau)) = - \lambda_N
\ee
with $\lambda_N$ the principal Dirichlet eigenvalue of $-\Delta/2$ on
$\Lambda_N$.
Combining (\ref{lb5})--(\ref{addlambda}), we get
\be{addP}
\lim_{\tau\ra\infty} \frac{1}{\tau} \log P(|W^{a\tau^{-1/(d-2)}}(\tau)| \leq b)
\geq - \widetilde I_N^{\kappa_a}(b) - \lambda_N \qquad \hbox{for all } N.\hskip.25in
\ee
Since $\lim_{N\ra\infty}\lambda_N=0$, it therefore suffices to show that
\be{lb7}
\lim_{N\ra\infty}\widetilde I_N^{\kappa_a}(b)=I^{\kappa_a}(b).
\ee
But this follows from the same type of argument as in Section 2.6.
\enddemo

\section{Upper and lower bounds in Theorem \ref{t2}}

In this section we explain how the arguments given in Sections 2 and 3 for the
Wiener sausage in $d\geq 3$ can be carried over to $d=2$. The necessary
modifications
are relatively minor and mainly involve a change in the choice of the
scaling parameters.

\demo{Upper bound} 
1.  Fix $N>0$. Wrap the Brownian motion around $\Lambda_{N\sqrt{t/\log
t}}$,
shrink space by $\sqrt{t/\log t}$ and time by $t/\log t$. Then the analogue of
(\ref{tau}), (\ref{ubscal}) reads
\be{ul1}
P(|W^a(t)| \leq bt/\log t) \leq P(|W_N^{a\sqrt{\tau e^{-\tau}}}(\tau)| \leq
b) \mbox{ with } \tau = \log t.\hskip.5in
\ee
We shall show how to obtain the analogue of Proposition \ref{p1}, namely,
\be{ul2}
\lim\limits_{\tau \ra \infty} \frac{1}{\tau} \log P(|W_N^{a\sqrt{\tau
e^{-\tau}}}(\tau)| \leq b)
= - I^{2\pi}_N(b), 
\ee
where $I^{2\pi}_N(b)$ is given by the same formulas as in
(\ref{mdsrf}), (\ref{mdsset}),
except that $\R^d$ is replaced by $\Lambda_N$ and $\kappa_a$ by $2\pi$. Since,
in Section 2.6, Proposition \ref{p2} was actually proved for any dimension, the
claim in (\ref{ul2}) will provide the upper bound in Theorem \ref{t2}.\\

 2.  Henceforth we suppress the indices $a,N$ and abbreviate
\be{pp1}
V_{\tau}=|W_N^{a\sqrt{\tau e^{-\tau}}}(\tau)|.
\ee
The analogue of Proposition \ref{VLDP} in Section 2.1 for $d=2$ reads:
\enddemo

\specialnumber{7} \proclaim{Proposition} \label{VLDPd=2}
$(V_{\tau})_{\tau>0}$ satisfies the {\rm LDP} on $\R_+$ with
rate $\tau$ and with rate function
\be{Jd=2}
J^{2\pi}_N(b)=\inf_{\phi \in \partial \Phi^{2\pi}_N(b)}\left[
\frac{1}{2} \int_{\Lambda_N} |\nabla\phi|^2(x) dx \right]
\ee
with
\be{Jsetd=2}
\partial \Phi^{2\pi}_N(b) = \left\{\phi\in H^1(\Lambda_N) \colon \,
\int_{\Lambda_N} \phi^2(x)dx = 1, \,\int_{\Lambda_N}
\Big(1-e^{-2\pi \phi^2(x)}\Big)dx = b \right\}.
\ee
\ep

\noindent This is the same as Proposition \ref{VLDP}, but with $\kappa_a$ replaced by
$2\pi$.
To prove Proposition \ref{VLDPd=2}, the coarse-graining argument in
Sections 2.2--2.4
can essentially be  copied. All that we need to do is replace $T_{\tau}$
defined in
(\ref{Tdef}) everywhere by
\be{Tdefd=2}
T_{\tau} = \frac{1}{\tau e^{-\tau}}
\ee
and prove the technical lemmas.
\medbreak

 3.  Section 2.2 carries over with the following difference.
 On the
right-hand side
of (\ref{Vhat3}) we end up with the expression
\be{Vhat3**}
E\Big(\exp\left[\frac{\tau}{\varepsilon T_\tau}
\left|W^a(\varepsilon T_\tau)\right|\right]\Big),
\ee
i.e., with an extra factor $\tau$ in the exponent. Since $\tau/\log
T_{\tau} \ra 1$ as
$\tau\ra\infty$, this means that instead of (\ref{boundheat}) we now need that
\be{boundheatd=2}
\sup_{T \geq 1} E\Big(\exp\left[\frac{\log T}{T}|W^a(T)|\right]\Big)<\infty.
\ee
However, this again follows from the results in van den Berg and Bolthausen
\cite{BeBo}.
Also, on the right-hand side of (\ref{10}) we end up with the expression
\be{mm1}
E_{x\sqrt{\varepsilon T_{\tau}},\varepsilon
T_{\tau}}\biggl(\exp\biggl[\frac{\tau}{\varepsilon T_{\tau}}
|W^a(\varepsilon T_{\tau})|\biggr]\biggr),
\ee
i.e., again with the extra factor $\tau$ in the exponent. This too can be
accommodated
because of (\ref{boundheatd=2}).
\medbreak

 4. Section 2.3 carries over after we prove Lemmas
\ref{cut1}--\ref{cont} for the new scaling
in (\ref{Tdefd=2}), with the following difference. We need to adapt the
argument\break at the
point where we are cutting out small holes around the points $\beta(i\varepsilon)$,
$1 \leq i \leq \tau/\varepsilon_{\phantom{\int}}\!\!$ (recall (\ref{18}), (\ref{cut})). Namely, this
time we cut
out holes of radius $1/\sqrt{\log T^{\phantom{1}}_\tau \log\log T_\tau}$, which is
considerably
larger than the radius $K/\sqrt{T_\tau}$ used before. The total volume of
the holes is at most $(\tau/\varepsilon+1)(\pi/\log T_\tau \log\log T_\tau)$, which for
$\tau \ra \infty$ tends to zero and therefore is negligible. The larger
radius is
needed for Part (a) of the new version of Lemma \ref{cut1}, which reads:

\specialnumber{8} \proclaim{Lemma} \label{cutd=2}
{\rm (a)} $\lim_{\tau \ra \infty} \sup_{y,z \notin B_{1/\sqrt{\log T_\tau
\log\log T_\tau}}}
q_{\tau,\varepsilon}(y,z)=0${\rm .}

{\rm (b)} $\lim_{\tau \ra \infty} \sup_{y,z \notin B_\rho}
|\tau q_{\tau,\varepsilon}(y,z)-2\pi\varphi_\varepsilon(y,z)| = 0$ for all
$0<\rho<N/4${\rm .}
\el

\demo{Proof}
(a) Step i of Part (a) in the proof of Lemma \ref{cut1} carries over, so that
(\ref{decom5})
again applies. Step ii of Part (a) is replaced by the following argument.
For any
$R>|y|>eb>0$,
\be{aii}
P_y^\infty(\sigma_b \leq \varepsilon/2) \leq P_y^\infty(\sigma_b <
\widehat\sigma_R)
+ P_y^\infty(\widehat\sigma_R \leq \varepsilon/2),
\ee
where $\widehat\sigma_R$ is the first entrance time into
$B^c_R=\R^d\setminus B_R$.
We have
\begin{eqnarray}\label{aii2}
P_y^\infty(\sigma_b < \widehat\sigma_R)
&=& \log \left(\frac{R}{|y|}\right) \Big/ \log \left(\frac{R}{b}\right), \\
P_y^\infty(\widehat\sigma_R \leq \varepsilon/2)
&\leq& 4 \exp\left[-\frac{(R-|y|)^2}{\varepsilon}\right]
\nn 
\end{eqnarray}
(for the latter see e.g.\ van den Berg and Davies \cite[Lemma 6.3]{vdBDa}).
The choice $R=|y|+\sqrt{\varepsilon \log\log (|y|/b)}$ together with the
inequality $\log (1+x) \leq x$ for $x \geq 0$ yields
\be{aii3}
P_y^\infty(\sigma_b \leq \varepsilon/2) \leq
\frac{1}{\log (\frac{|y|}{b})} \left[4 + \frac{1}{|y|}
\sqrt{\varepsilon \log\log \left(\frac{|y|}{b}\right)}\ \right].
\ee
Inserting this into (\ref{decom5}), we get for any $0<eb<b'<N/4$,
\be{aii4}
\sup_{y,z \notin B_{b'}} q_b(y,z) \leq \frac{2c_2}{1-c_3}
\frac{1}{\log \left(\frac{b'}{b}\right)} \left[4 + \frac{1}{b'}
\sqrt{\varepsilon \log\log \left(\frac{b'}{b}\right)}\ \right].\hskip.5in
\ee
Now put $b=a/\sqrt{T_\tau}$, $b'=1/\sqrt{\log T_\tau \log\log T_\tau}$ and
use that $\log T_\tau \sim \tau$\break ($\tau \to \infty$), to get the claim.\medbreak

(b) Part (b) in the proof of Lemma \ref{cut1} carries over, with the only
difference that (\ref{LG1}) is to be replaced by
\be{LG2}
\lim_{b \da 0} \frac{1}{\pi/\log(\frac{1}{b})} P_y^\infty(\sigma_b \leq t)
= \int_0^t p_s(-y) ds \qquad \mbox{ for all } y \in \R^2, t \geq 0\hskip.25in
\ee
(Le Gall \cite{LG86}). For $b=a/\sqrt{T_\tau}$ we have $\pi/\log(\frac{1}{b})
\sim 2\pi/\tau$ ($\tau \to \infty$), which explains how the factor $2\pi$
arises, replacing $\kappa_a$.
\enddemo

Lemmas \ref{cont1} and \ref{cont} were based on Lemma \ref{cut1}(b). It is obvious
that with the new version in Lemma \ref{cutd=2}(b) the rest of the argument in
Section 2.3 is unchanged.
\medbreak

  5.  Section 2.4 carries over verbatim with only $\kappa_a$ to be
replaced by $2\pi$
everywhere. Section 2.5 also has no changes. In fact, in both these
sections dimension
plays no role at all.

\demo{Lower bound} 
The proof of Proposition \ref{plb} carries over after the appropriate changes in
  scaling. \enddemo

\section{Analysis of the variational problem}

This section contains the main analytic part of our paper. Theorems
\ref{t3}(i)--(iii) are proved
in Sections 5.1--5.3, Theorems \ref{t4}(i), (ii) in Sections 5.4, 5.5, and
Theorems \ref{t5}(i)--(iii)
in Sections 5.6--5.8. Recall the notation introduced in Section~1.

We will repeatedly make use of the following scaling relations. Let $\phi
\in H^1(\R^d)$.
For $p,q>0$, define $\psi \in H^1(\R^d)$ by
\be{scal*}
\psi(x) = q\phi(x/p).
\ee
Then
\begin{eqnarray}\label{scal**}
&&\|\nabla\psi\|^2_2 = q^2p^{d-2} \|\nabla\phi\|^2_2, \quad
\|\psi\|^2_2 = q^2p^d \|\phi\|^2_2,\quad
\|\psi\|^4_4 = q^4p^d \|\phi\|^4_4,  \\[4pt]
&&\hskip1in\int (1-e^{-\psi^2}) = p^d \int (1-e^{-q^2 \phi^2}).
\nn
\end{eqnarray}

We will also repeatedly make use of the following Sobolev inequalities (see
Lieb and Loss \cite[pp.\ 186 and 190]{LL}):
\be{4.38}
S_d\|f\|_q^2 \leq \|\nabla f\|_2^2 \qquad (d \geq 3, f\in D^1(\R^d) \cap
L^2(\R^d))
\ee
with
\be{4.40}
q=\frac{2d}{d-2},\quad
S_d=d(d-2)2^{-2(d-1)/d}\pi^{(d+1)/d}\Big[\Gamma\Big(\frac{d+1}{2}\Big)\Big]^
{-2/d},\hskip.35in
\ee
and
\be{4.36}
\|f\|_4 \le S_{2,4}(\|\nabla f\|_2^2+\|f\|_2^2)^{1/2} \qquad (d=2, f \in
H^1(\R^2))
\ee
with $S_{2,4} = (4/27\pi)^{1/4}$.

\demo{{\rm 5.1.} Proof of Theorem {\rm \ref{t3}(i),} reduction to radially symmetric nonincreasing
functions{\rm ,} and adaptation of the constraints}
\medbreak
 1.  We begin by reformulating the variational problem for
$I^{\kappa_a}(b)$ in Theorem
\ref{t1}.

\specialnumber{9} \proclaim{Lemma} \label{l1}
Let $d\ge2$ and $a>0${\rm .} For every $b>0$
\be{4.1}
I^{\kappa_a}(b)= \frac{1}{2 \kappa_a^{2/d}}
\chi\biggl(\frac{b}{\kappa_a}\biggr),
\ee
where $\chi \colon (0,\infty) \mapsto [0,\infty)$ is given by
\be{4.2}
\chi(u)=\inf\biggl\{\|\nabla\psi\|_2^2 \colon ~\psi \in
H^1(\R^d),\,\|\psi\|_2=1,
\,\int_{\R^d} (1-e^{-\psi^2})\le u\biggr\}.\hskip.35in
\ee
\el

\demo{Proof}
Apply (\ref{scal*}) and (\ref{scal**}) with $p=\kappa_a^{-1/d}$ and
$q=\kappa_a^{1/2}$
to (\ref{mdsrf}), (\ref{mdsset}).
\enddemo

 \noindent Lemma \ref{l1} proves Theorem \ref{t3}(i).

\medbreak 2.  The following lemma reduces the variational problem in (\ref{4.2})
to radially
symmetric nonincreasing ($\RSNI$) functions. This reduction will become
important later on.

\specialnumber{10} \proclaim{Lemma} \label{l2}
Let
\be{4.6}
\cR_u=\biggl\{\psi \in H^1(\R^d) \colon ~\psi ~\RSNI,\,
\|\psi\|_2=1,\,\int_{\R^d}(1-e^{-\psi^2})\le u\biggr\}.\hskip.25in
\ee
Then
\be{4.7}
\chi(u)=\inf\{\,\|\nabla\psi\|_2^2 \colon ~\psi\in\cR_u\,\}.
\ee
\el

\demo{Proof}
It is clear that $\chi(u)\le\inf\{\|\nabla\psi\|_2^2 \colon ~\psi\in \cR_u\}$.
To prove the reverse, we let $\psi^*$ denote the symmetric decreasing
rearrangement
of $\psi$. Then (see Lieb and Loss \cite[\S 3.3 and \S 7.17]{LL})
$\psi^*$ is
nonnegative, $\RSNI$, and \vglue-9pt
\be{4.9}
\|\nabla\psi\|_2\ge\|\nabla\psi^*\|_2, \qquad
\|\psi\|_2=\|\psi^*\|_2, \qquad
\int(1-e^{-\psi^2})=\int(1-e^{-\psi^{*2}}).
\ee
Hence
\begin{eqnarray} \label{4.10}
&&\\
\chi(u) &\geq& \inf\{\,\|\nabla\psi^*\|_2^2 \colon ~\psi \in H^1(\R^d),
~\|\psi\|_2=1,
\int(1-e^{-\psi^2})\le u\,\}\nn  \\[4pt]
&=& \inf\{\,\|\nabla\psi^*\|_2^2 \colon ~\psi^* \in H^1(\R^d), ~\|\psi^*\|_2=1,
\int(1-e^{-\psi^{*2}})\le u\,\}\nn \\[4pt]
&\geq& \inf\{\,\|\nabla\psi\|_2^2 \colon ~\psi\in \cR_u\,\}.
\nn
\end{eqnarray}
\vglue-24pt
\enddemo

\phantom{away}

  3.  The following lemma makes a statement about the minimisers of
(\ref{4.2}).
Whether or not these exist will be established later on.

\specialnumber{11} \proclaim{Lemma} \label{RSNI}
Any minimiser of {\rm (\ref{4.2})} is strictly positive{\rm ,} radially symmetric
{\rm (}\/modulo shifts\/{\rm )}
and strictly decreasing in the radial component{\rm .}
\el

\phantom{line away}

\demo{Proof}
Let $\psi$ be any minimiser of (\ref{4.2}). Let $\psi^*$ be its symmetric
decreasing
rearrangement. Then, by (\ref{4.9}), $\psi^*$ too is a minimiser of
(\ref{4.2}).
By Brothers and Ziemer \cite[Th.~1.1]{BZ}, $\|\nabla\psi\|_2 >
\|\nabla\psi^*\|_2$
if $\psi$ is not a shift of $\psi^*$ and the set $\{x \in \R^d
\colon\,(\nabla\psi^*)(x)=0\}$
has zero Lebesgue measure. We will show that $d\psi^*/dr<0$. Therefore
$\psi$ must be
a shift of $\psi^*$ (otherwise $\psi$ could not be a minimiser) and the
claim will follow.

Since $\psi^*$ is a radially symmetric minimiser of (\ref{4.2}), it
satisfies the
Euler-Lagrange equation \vglue-9pt
\be{ELRSNI}
\frac{d^2\psi^*}{dr^2} + \frac{d-1}{r}\frac{d\psi^*}{dr} =
\lambda \psi^*(1-e^{-\psi^{*2}}) + \mu \psi^* \quad (r>0),
\ee
where $\lambda,\mu$ are Lagrange multipliers (see Berestycki and Lions
\cite[\S 5b]{BeLi}). By differentiating (\ref{ELRSNI}) repeatedly with respect to $r$,
we see that $\psi^* \in C^{\infty}(0,\infty)$. Now, we already know that
$d\psi^*/dr\leq 0$. Suppose that $(d\psi^*/dr)(r_0)=0$ for some $r_0>0$.
Then clearly
we must also have $(d^2\psi^*/dr^2)(r_0)=0$. But from the derivatives of
(\ref{ELRSNI})
it then  follows that $(d^n\psi^*/dr^n)(r_0)=0$ for all $n\in\N$. However,
(\ref{ELRSNI})
is a second order differential equation with Lipschitz coefficients, and
therefore the
latter entails that $\psi^*(r)=\psi^*(r_0)$ for all $r\in(0,\infty)$, i.e.,
$\psi^*$ is
constant. But this contradicts $\|\psi^*\|_2=1$. Hence
$(d\psi^*/dr)(r_0)<0$. This proves
the claim since $r_0$ was arbitrary.
\enddemo

  4.  We end this section with a lemma stating that the constraints in
(\ref{4.2})
can be adapted. This will turn out to be important later on.

\specialnumber{12} \proclaim{Lemma} \label{l6}
Let
\begin{eqnarray}\label{eq4.20}
\widehat \chi(u) &=& \inf\{\|\nabla
\psi\|_2^2 \colon ~\|\psi\|_2=1,\,\int_{\R^d}(1-e^{-\psi^2})=u\},  \\[4pt]
\widetilde \chi(u) &=& \inf\{\|\nabla
\psi\|_2^2 \colon ~\|\psi\|_2\leq1,\,\int_{\R^d}(1-e^{-\psi^2})=u\}.
\nn
\end{eqnarray}
Then
\be{eq4.22}
\chi(u)= \widehat \chi(u) = \widetilde \chi(u).
\ee
\el

\demo{Proof}
We use an approximation argument.

\medbreak i. It is clear from (\ref{4.2}) and (\ref{eq4.20}) that
$\chi(u)\le \widehat \chi(u)$. To prove the reverse, let $(\psi_j)$ be a
minimising
sequence of $\chi(u)$. Then $\|\psi_j\|_2=1$, $\int(1-e^{-\psi_j^2})\le u$, and
$\|\nabla\psi_j\|^2_2 \to \chi(u)$ as $j \to \infty$. Define, for $a>0$,
\be{eq4.26}
g_{\psi}(a)=a^d\int(1-e^{-a^{-d}\psi^2}).
\ee
Then
\be{eq4.27}
g_{\psi}'(a)=da^{d-1}\int(1-e^{-a^{-d}\psi^2}-a^{-d}\psi^2
e^{-a^{-d}\psi^2}).
\ee
Since $1-e^{-x}-xe^{-x}\ge0$ for $x\geq 0$, we have that
$g_{\psi}'(a)\ge0$. Since
$g_{\psi_j}(\infty)=\|\psi_j\|_2^2=1$ and $g_{\psi_j}(1)=
\int(1-e^{-\psi_j^2}) \leq u$,
we see that there exists a sequence $(a_j)$ with $a_j\ge 1$ such that
\be{eq4.28}
g_{\psi_j}(a_j) = u \quad \mbox{ for all } j.
\ee
Next, let $\phi_j\in H^1(\R^d)$ be defined by
$\phi_j(x)=a_j^{-d/2}\psi_j(x/a_j)$.
Then, recalling (\ref{scal*}), (\ref{scal**}) and using (\ref{eq4.28}), we
see that \vglue-9pt
\be{Phiprop}
\|\nabla\phi_j\|^2_2 = \frac{1}{a_j^2}\|\nabla \psi_j\|^2_2,\quad
\|\phi_j\|^2_2=\|\psi_j\|^2_2=1,\quad
\int(1-e^{-\phi_j^2}) = u \quad \mbox{for all } j.
\ee
Hence (\ref{eq4.20}) gives 
\be{eq4.31}
\widehat \chi (u)\le\|\nabla\phi_j\|_2^2=\frac1{a_j^2}\|\nabla
\psi_j\|_2^2\le\|\nabla\psi_j\|_2^2 \quad \mbox{for all } j.
\ee
But $\|\nabla\psi_j\|_2^2\to\chi(u)$ as $j\to\infty$, and so $\widehat
\chi(u)\le\chi(u)$.\smallbreak

  ii.  It is clear from (\ref{eq4.20}) that $\widetilde \chi(u)\le
\widehat \chi(u)$. To
prove the reverse, we begin with the following observation:

\specialnumber{13} \proclaim{Lemma} \label{lcompset}
The set
\be{compset}
\Big\{\psi \in H^1(\R^d) \colon ~\psi \,\RSNI,\,\|\nabla\psi\|_2 \leq C,
\,\|\psi\|_2 \leq 1,\,\int(e^{-\psi^2}-1+\psi^2) = 1-u \Big\}
\ee
is a compact for all $C<\infty$.
\el

Before proving Lemma \ref{lcompset} we first complete the argument. Since
$\psi \mapsto
\|\nabla\psi\|_2$ is lower semi-continuous, it follows from (\ref{compset})
that the
variational problem for $\widetilde \chi(u)$ has a minimiser, say $\psi^*$.
For $n \in \N$, let
\be{pdef1}
p_n(x) = \frac{1}{\pi^{d/2}n^d}e^{-|x|^2/n^2} \quad (x \in \R^d),
\ee
and note that $\int p_n = 1$ and $\int (\nabla\sqrt{p_n})^2 = 2d/n^2$. Now
define $\psi^*_n$ by
\be{psidef1}
\psi^{*2}_n = \psi^{*2} + [1-\|\psi^*\|^2_2] ~p_n.
\ee
Then $\|\psi^*_n\|_2=1$ for all $n$. Moreover, since $x \mapsto e^{-x}-1+x$
is increasing on
$[0,\infty)$, we have
\be{intext1}
\int (e^{-\psi^{*2}_n}-1+\psi^{*2}_n) \geq
\int (e^{-\psi^{*2}}-1+\psi^{*2}) = 1-u \quad \mbox{for all } n.\hskip.5in
\ee
Therefore $\psi^*_n$ satisfies the constraints in the variational problem
for $\chi(u)$,
implying that
\be{chiest1}
\chi(u) \leq \|\nabla\psi^*_n\|^2_2 \quad \mbox{ for all } n.
\ee
But by the convexity inequality for gradients (Lieb and Loss \cite[Th.\ 7.8]{LL}) we have
\be{gradest1}
\|\nabla\psi^*_n\|^2_2 \leq \|\nabla\psi^*\|^2_2 + [1-\|\psi^*\|^2_2]\,
\|\nabla\sqrt{p_n}\|^2_2 = \widetilde \chi(u) + [1-\|\psi^*\|^2_2]\,
\frac{2d}{n^2}.\enspace
\ee
Letting $n \ra \infty$, we thus end up with $\chi(u)\leq\widetilde\chi(u)$. But
$\chi(u)=\widehat\chi(u)$ by Step i, and so the claim is proved.\medbreak

iii.  It thus remains to prove Lemma
\ref{lcompset}.

\demo{Proof}
The key point is to show that the contribution to the integral in
(\ref{compset})
coming from small $x$ and from large $x$ is uniformly small. First we pick
$0<R<\infty$ and estimate
\be{Rest}
\int_{B^c_R} (e^{-\psi^2}-1+\psi^2) \leq \int_{B^c_R} \frac{1}{2}\psi^4
\leq \frac{1}{2\omega_d R^d} \int_{B^c_R} \psi^2 \leq \frac{1}{2\omega_d R^d},\hskip.5in
\ee
where $\omega_d$ is the volume of the ball with unit radius, and we use that
$\psi$ is $\RSNI$ and $\|\psi\|_2 \leq 1$. Next, we pick $0<r<\infty$ and
estimate,
using H\"older's inequality,
\be{rest}
\int_{B_r} (e^{-\psi^2}-1+\psi^2) \leq \int_{B_r} \psi^2
\leq (\omega_d r^d)^{1/p} \|\psi\|^2_{2q}  \qquad \Big(p,q \geq
1,\frac{1}{p}+\frac{1}{q}
= 1 \Big).
\ee
The last factor may be estimated with the help of the Sobolev inequalities in
(\ref{4.38})--(\ref{4.36}):
\begin{eqnarray}\label{Sobolev}
\|\psi\|^2_{2d/(d-2)} \enspace\ \leq  &S_d \|\nabla\psi\|^2_2\phantom{whatever1}  &(d \geq 3)\\[4pt]
\|\psi\|^2_4 \enspace\ \leq & S_{2,4} (\|\nabla\psi\|^2_2 + \|\psi\|^2_2) &(d=2).
\nn
\end{eqnarray}
Thus, picking $p=d/2,q=d/(d-2)$ for $d\geq 3$ and $p=q=2$ for $d=2$, we obtain
using $\|\nabla\psi\|_2 \leq C$ that
\be{rest2}
\int_{B_r} (e^{-\psi^2}-1+\psi^2) \leq
\left\{\ba{ll}
C_d r        &(d \geq 3)\\[4pt]
C_2 r^{1/2}  &(d=2).
\ea
\right.
\ee
We see from (\ref{Rest}) and (\ref{rest2}) that the contribution from
$B^c_R$ and
$B_r$ tends to zero uniformly in $\psi$ as $R \ra \infty$ and $r \da 0$. We can
now complete the proof as follows. Any sequence $(\psi_j)$ in $H^1(\R^d)$ has a
subsequence that converges to some $\psi \in H^1(\R^d)$ uniformly on every
annulus
$B_R \setminus B_r$ (use the fact that $\psi_j$ is $\RSNI$ and $\|\psi_j\|_2\leq 1$
for all
$j$). Because $\psi_j$ is $\RSNI$, $\|\nabla\psi_j\|_2\leq C$,
$\|\psi_j\|_2\leq 1$
for all $j$, the same is true for $\psi$. Moreover, since
\begin{eqnarray}\label{psiint}
&&\int (e^{-\psi_j^2}-1+\psi_j^2) = 1-u \quad \mbox{ for all } j,  \\
&&\lim\limits_{j \ra \infty} \int_{B_R \setminus B_r} (e^{-\psi_j^2}-1+\psi_j^2)
= \int_{B_R \setminus B_r} (e^{-\psi^2}-1+\psi^2),\nn \\
&&\lim\limits_{r \da 0, R \ra \infty} \int_{B_R \setminus B_r}
(e^{-\psi^2}-1+\psi^2)
= \int (e^{-\psi^2}-1+\psi^2),
\nn
\end{eqnarray}
we also have $\int (e^{-\psi^2}-1+\psi^2)=1-u$. Therefore $\psi$ is in the set.
\enddemo

This completes the proof of Lemma \ref{lcompset} and hence of Lemma \ref{l6}.\hfill\qed
\enddemo

The reason behind Lemma \ref{lcompset} is the following. Although $\psi$
may lose
$L^2$-mass to infinity, the integral cannot. Indeed, following an argument in
Brezis and Lieb \cite{BL}, we can show that if $\|\psi\|_2^2$ loses mass $\rho
\in (0,u]$, then also $\int (1-e^{-\psi^2})$ loses mass $\rho$, and so $\int
(e^{-\psi^2} - 1 + \psi^2)$ loses nothing.

In the sequel we shall often suppress the condition $\psi \in H^1(\R^d)$ from
the notation.

\demo{{\rm 5.2.} Proof of Theorem {\rm \ref{t3}(ii)}}
\enddemo

 1.  Since $1-e^{-x}\le x$ for $x \ge 0$, we have $\int (1-e^{-\psi^2})
\le \|\psi\|_2^2$.
So, for $u\ge 1$, (\ref{4.2}) reduces to
\be{4.12}
\chi(u)=\inf\{\,\|\nabla\psi\|_2^2 \colon ~\|\psi\|_2=1\,\}.
\ee
Suppose $\psi \in H^1(\R^d)$ is such that $\|\psi\|_2=1$. Apply to
(\ref{4.12}) the scaling
in (\ref{scal*}) and (\ref{scal**}) with $p>0$ arbitrary and $q=p^{-d/2}$, to
obtain $\chi(u)\le p^{-2}
\|\nabla\psi\|_2^2$ for $u \ge 1$. Taking the limit $p\to\infty$, we get
$\chi(u)=0$ for $u\geq 1$.\medbreak

 2.  It follows from Theorems \ref{t4}(ii) and \ref{t5}(iii) that
$\chi$ is strictly positive
in a left-neighbourhood of 0. Since, by Theorem \ref{t3}(iii), $u \mapsto
u^{2/d}\chi(u)$ is
nonincreasing on $(0,1)$, it follows that $\chi$ is strictly decreasing on
$(0,1)$.\medbreak

  3.  Step 1 shows that $\chi$ is continuous on $(1,\infty)$. Theorems
\ref{t4}(ii) and
\ref{t5}(iii) imply that $\chi$ is continuous at $u=1$. Therefore we need
only prove continuity
on $(0,1)$. Let $u_0\in(0,1)$ be arbitrary. Since $\chi$ is lower
semi-continuous and
nonincreasing, it is right-continuous. Let
\be{eq4.39}
\delta=\lim_{u\ua u_0}\chi(u)-\chi(u_0) \geq 0.
\ee
We shall show that $\delta = 0$ by using a perturbation argument.\medbreak

4.  Let $\varepsilon>0$ be arbitrary. Then,
because
$\chi(u) = \widehat
\chi(u)$
by Lemma \ref{l6}, there exist $\psi_\varepsilon,\Phi_\varepsilon\in H^1(\R^d)$
satisfying
\bel{eq4.40}
\quad \|\psi_\varepsilon\|_2=1, \int(1-e^{-\psi_\varepsilon^2})=u_0,\phantom{\|\psi_1}
&&\|\nabla\psi_\varepsilon\|_2^2 \le\chi(u_0)+\varepsilon,  \\
\quad \|\Phi_\varepsilon\|_2=1, \int(1-e^{-\Phi_\varepsilon^2}) = u_0-\varepsilon,
&&\|\nabla\Phi_\varepsilon\|_2^2 \le\chi(u_0-\varepsilon)+\varepsilon.
\nne
Define, for $0\le\alpha\le1$,
\be{eq4.42}
\Lambda_{\alpha,\varepsilon}=[\alpha\psi_\varepsilon^2
+(1-\alpha)\Phi_\varepsilon^2]^{1/2}.
\ee
Then, by (\ref{eq4.40}),
\be{eq4.43}
\|\Lambda_{\alpha,\varepsilon}\|_2^2=
\int[\alpha\psi_\varepsilon^2+(1-\alpha)\Phi_\varepsilon^2]=1
\ee
and, by the convexity inequality for gradients (see Lieb and Loss \cite[Th.\ 7.8]{LL}), 
\begin{eqnarray}\label{eq4.44}
\|\nabla\Lambda_{\alpha,\varepsilon}\|_2^2
&\le&\alpha\|\nabla\psi_\varepsilon\|_2^2+(1-\alpha)\|\nabla
\Phi_\varepsilon\|_2^2  \\[4pt]
&\le&\alpha(\chi(u_0)+\varepsilon)+(1-\alpha)(\chi(u_0-\varepsilon)+\varepsilon)\nn \\[4pt]
&=&\alpha\chi(u_0)+(1-\alpha)\chi(u_0-\varepsilon)+\varepsilon.
\nn
\end{eqnarray}
Next define, for $0\le\alpha\le1$,
\be{eq4.45}
k(\alpha)=\int(1-e^{-\Lambda^2_{\alpha,\varepsilon}}).
\ee
Then, by (\ref{eq4.42}),
\be{eq4.46}
k''(\alpha)=\int(\psi_\varepsilon^2-\Phi_\varepsilon^2)^2
e^{-\alpha\psi_\varepsilon^2-(1-\alpha)\Phi_\varepsilon^2}.
\ee
It follows that $k$ is convex on $[0,1]$, and consequently
\be{eq4.47}
k(\alpha) \le \alpha k(1) + (1-\alpha) k(0) = \alpha u_0 +
(1-\alpha)(u_0-\varepsilon).
\ee
By the convexity of $k$ and by (\ref{eq4.47}), there exists a unique
$\alpha_\varepsilon\in[\frac12,1)$ such that $k(\alpha_\varepsilon)=u_0-\varepsilon/2$.
By (\ref{eq4.43})--(\ref{eq4.45}), Lemma \ref{l6} and the fact that $\chi$
is nonincreasing, we therefore have
\begin{eqnarray}\label{eq4.49}
\qquad\ \ \chi(u_0-\varepsilon/2) &  = & \widehat \chi(u_0-\varepsilon/2)
\le\|\nabla \Lambda_{\alpha_\varepsilon,\varepsilon}\|_2^2\\
&  \le &\alpha_\varepsilon\chi(u_0) 
+\ (1-\alpha_\varepsilon)
\chi(u_0-\varepsilon)+\varepsilon\nn\\
&\le&\displaystyle\frac12\chi(u_0)+\displaystyle\frac12
\chi(u_0-\varepsilon)+\varepsilon.
\nn\end{eqnarray}
Hence
\be{eq4.50}
\lim_{\varepsilon \da 0}\chi(u_0-\varepsilon/2)\le\frac12\chi(u_0)
+\frac12\lim_{\varepsilon \da 0}\chi(u_0-\varepsilon).
\ee
Combining (\ref{eq4.39}) and (\ref{eq4.50}), we therefore arrive at
\be{eq4.51}
\chi(u_0)+\delta\le\frac12\chi(u_0)+\frac12[\chi(u_0)+\delta].
\ee
Thus $\delta=0$ and $\chi$ is continuous at $u_0$. This proves the
continuity of $\chi$ on $(0,1)$, since $u_0\in(0,1)$ was arbitrary.

\demo{{\rm 5.3.} Proof of Theorem {\rm \ref{t3}(iii)}}
\medbreak  
1.  The first claim in Theorem \ref{t3}(iii) is proved as follows.
Apply (\ref{scal*}), (\ref{scal**}) with $p=u^{1/d}$ and $q=u^{-1/2}$
to (\ref{4.2}), to obtain
\be{chiscal3}
u^{2/d}\chi(u) = \inf\{\|\nabla \psi\|^2_2 \colon ~\|\psi\|_2 = 1,
~\textstyle{\int} (1-e^{-u^{-1}\psi^2}) \leq 1 \}.\hskip.5in
\ee
Since the integrand is nonincreasing in $u$, so is the infimum. Therefore
we find that $u \mapsto u^{2/d}\chi(u)$ is nonincreasing. To prove the strict
monotonicity claimed in Theorem \ref{t3}(iii), we need to wait until the end
of Section 5.8, as this will require the existence of a minimiser for a certain
range of $u$--values. \medbreak

2.  The second claim in Theorem \ref{t3}(iii) is proved by deriving
upper and lower bounds. Pick $\delta>0$. Let $B$ be any open set in $\R^d$
with $|B|=1$. Then, since $1-e^{-\delta^{-1}\psi^2}\ge 0$, it follows from
(\ref{chiscal3}) that
\bel{4.17}&&\\
\delta^{2/d}\chi(\delta) &\leq& \inf\biggl\{\|\nabla\psi\|_2^2 \colon
~\|\psi\|_2=1,\,
\int_{\R^d}(1-e^{-\delta^{-1}\psi^2})\le 1, \hbox{supp}(\psi)\subset B
\biggr\}  \nn\\
&=& \inf\biggl\{\|\nabla\psi\|_2^2 \colon ~\|\psi\|_2=1,\,
\int_B(1-e^{-\delta^{-1}\psi^2})\leq 1,\,\hbox{supp}(\psi)\subset B
\biggr\}\nn \\
&\leq& \inf\{\,\|\nabla\psi\|_2^2 \colon
~\|\psi\|_2=1,\,\hbox{supp}(\psi)\subset B \,\}.
\nne
Take the infimum over $B$ of the right-hand side of (\ref{4.17}). This
infimum equals
$\lambda_d$ by the Faber-Krahn isoperimetric inequality for the Dirichlet
Laplacian
(see Faris \cite{F}). Hence $\delta^{2/d} \chi(\delta)\leq \lambda_d$,
which proves the upper bound.\medbreak

 3.  The lower bound is more laborious. Let $\psi$ be a minimiser of
(\ref{chiscal3}).
(In Sections 5.4 and 5.7 we will prove the existence of a minimiser of
(\ref{4.2}), and hence
of (\ref{chiscal3}), for $u$ in a neighbourhood of 0 for any $d \ge 1$.) We
exploit the
radial symmetry and monotonicity of $\psi$ established in Lemma \ref{RSNI}.
For $t>0$,
define the ball $B_t=\{\,x\in\R^d \colon ~\psi^2(x)>t\,\}$ and put
$\mu(t)=|B_t|$.
Then
\be{4.21}
1 = \int_{\R^d} \psi^2 = \int_{[0,\infty)}t\,d[-\mu(t)].
\ee
Moreover, abbreviating $\varepsilon=\sqrt{\delta}$ we have
\bel{4.20}&&\\
1 \geq \int_{\R^d}(1-e^{-\delta^{-1}\psi^2})
&=&\int_{[0,\infty)}(1-e^{-\delta^{-1}t})\,d[-\mu(t)]\nn\\
&\geq &\int_{[\varepsilon,\infty)}(1-e^{-\delta^{-1}t})\,d[-\mu(t)]\nn  \\
&\geq& (1-e^{-1/\varepsilon})\int_{[\varepsilon,\infty)}\,d[-\mu(t)]
=(1-e^{-1/\varepsilon})|B_{\varepsilon}|.
\nne
Hence
\bel{4.22}
\quad\int_{\R^d\setminus B_{\varepsilon}}\psi^2
&= &\int_{[0,\varepsilon)}t\,d[-\mu(t)]
\\
&\le& \frac{\varepsilon}{1-e^{-1/\varepsilon}}
\int_{[0,\varepsilon)} (1-e^{-\delta^{-1}t})\,d[-\mu(t)]\le \frac{\varepsilon}{1-e^{-1/\varepsilon}},
\nne
where we use the first inequality in (\ref{4.20}). Combining (\ref{4.21})
and (\ref{4.22}) we obtain
\be{4.23}
\int_{B_{\varepsilon}}\psi^2 \ge1-\frac{\varepsilon}{1-e^{-1/\varepsilon}}.
\ee
Next, define $\zeta=\psi-\sqrt{\varepsilon}$. Then $\zeta>0$ on $B_{\varepsilon}$ and
$\zeta=0$ on $\partial B_{\varepsilon}$. By Cauchy-Schwarz and $\|\psi\|_2=1$,
we have
\be{4.25}
\int_{B_{\varepsilon}}\zeta\le\int_{B_{\varepsilon}}\psi\le\biggl(\int_{B_{\varepsilon}}
\psi^2\biggr)^{1/2}~|B_{\varepsilon}|^{1/2}\le|B_{\varepsilon}|^{1/2}.
\ee
Hence
\be{4.26}
\int_{B_{\varepsilon}} \psi^2 = \int_{B_{\varepsilon}}[\zeta+\sqrt{\varepsilon}]^2
\le \int_{B_{\varepsilon}}\zeta^2 + 2\sqrt{\varepsilon} |B_{\varepsilon}|^{1/2}
+\varepsilon |B_{\varepsilon}|.
\ee
Combining (\ref{4.20}), (\ref{4.23}), (\ref{4.26}) and
abbreviating $\eta = \varepsilon/(1-e^{-1/\varepsilon})$, we obtain
\be{4.28}
\int_{B_{\varepsilon}}\zeta^2 \ge 1 - 2 \sqrt{\eta} - 2 \eta.
\ee
{}Finally, if we define $\phi$ by $\zeta = (1-2\sqrt{\eta} -2\eta)^{1/2}\phi$,
then $\int_{B_{\varepsilon}} \phi^2\ge 1$ and $\phi\big|_{\partial
B_{\varepsilon}}=0$. So,
recalling that $\psi$ is a minimiser of (\ref{chiscal3}), we get\vglue-9pt
\be{4.31}
\delta^{2/d} \chi(\delta) = \int_{\R^d} |\nabla \psi|^2
\geq \int_{B_{\varepsilon}}|\nabla \psi|^2
= \int_{B_{\varepsilon}}|\nabla \zeta|^2
= (1-2\sqrt{\eta}-2\eta) \int_{B_{\varepsilon}}|\nabla \phi|^2
\ee
with
\bel{4.31*}
\int_{B_{\varepsilon}}|\nabla \phi|^2
&\geq& \inf\left\{\int_{B_{\varepsilon}}|\nabla\phi|^2 \colon
~\int_{B_{\varepsilon}}\phi^2\ge 1,\,\phi\big|_{\partial B_{\varepsilon}}=0\right\}  \\[4pt]
&=&\lambda_d(B_{\varepsilon}) = |B_{\varepsilon}|^{-2/d} \lambda_d
\ge (1-e^{-1/\varepsilon})^{2/d}\lambda_d,
\nne
where we use the scaling of the smallest Dirichlet eigenvalue of $-\Delta$ on
$B_{\varepsilon}$, in combination with (\ref{4.20}). Letting $\delta \da 0$
in (\ref{4.31}) and using the fact that $\varepsilon \da 0$ and $\eta \da 0$, we arrive at
$\liminf_{\delta\da 0}\delta^{2/d}\chi(\delta)\geq\lambda_d$.\enddemo

\demo{{\rm 5.4.} Proof of Theorem {\rm \ref{t4}(i)}}
To prove that for $2 \leq d \leq 4$ the variational problem
in (\ref{4.2}) has a minimiser for all $u \in (0,1)$, we do a variational
argument that takes advantage of Lemma \ref{l6}. Fix $u \in (0,1)$. Let
$\psi^*$ be any minimiser of $\widetilde \chi(u)$, which we know exists by
(\ref{compset}), i.e.,
\be{condext1}
\|\nabla\psi^*\|^2_2 = \widetilde \chi(u),\quad
\|\psi^*\|^2_2 \leq 1,\quad
\textstyle\int(e^{-\psi^{*2}}-1+\psi^{*2}) = 1-u.\hskip.5in
\ee
There are two cases. Either $\|\psi^*\|_2=1$, in which case $\psi^*$ is also
a minimiser of $\chi(u)$ and we are done, or $\|\psi^*\|_2<1$. It remains to
exclude the latter case.

\medbreak
\underline{$d=2$}: Suppose that $\|\psi^*\|_2^2=1-\rho$ ($\rho>0$).
Let $\phi^*(\cdot)=\psi^*(~\cdot~(1-\rho)^{1/2})$. Then, by
(\ref{scal*}), (\ref{scal**}), 
\be{eqsd=2}
\|\nabla\phi^*\|^2_2 = \|\nabla\psi^*\|^2_2 = \chi(u),\quad
\|\phi^*\|^2_2 = 1,\quad
\int (e^{-\phi^{*2}}-1+\phi^{*2}) = \frac{1-u}{1-\rho}.
\ee
Hence
\be{eqsd=2'}
\chi(u) = \|\nabla\phi^*\|^2_2 \geq \chi\Big(\frac{u-\rho}{1-\rho}\Big).
\ee
But $(u-\rho)/(1-\rho)<u$, and so the right-hand side is strictly larger than
$\chi(u)$ by Theorem \ref{t3}(ii), which is a contradiction.

\medbreak
\underline{$d=3,4$}:
We can do smooth perturbations of $\psi^*$ inside the class $$\{\psi \in
H^1(\R^d) \colon ~\|\psi\|_2 \leq 1,\, \int (e^{-\psi^2}-1+\psi^2)=1-u\}$$ to
conclude that $\psi^*$ must satisfy the Euler-Lagrange equation associated with
the variational problem
\be{varprl}
\inf\{\|\nabla\psi\|^2_2 \colon ~\textstyle\int (e^{-\psi^2}-1+\psi^2) = 1-u\}.
\ee
But then we have a contradiction to the following:

\specialnumber{14} \proclaim{Lemma} \label{EL1-4}
For $d=3,4$ all solutions of the Euler\/{\rm -}\/Lagrange equation associated with the
variational problem
\be{ELuc}
\inf\{\|\nabla\psi\|^2_2 \colon ~\textstyle\int (e^{-\psi^2}-1+\psi^2) = 1 \}
\ee
have infinite $L^2$\/{\rm -}\/norm{\rm .}
\el

\demo{Proof}
By the results of Section 5b in Berestycki and Lions \cite{BeLi}, there
exists a Lagrange multiplier $\lambda_d>0$ such that
\be{fax1**}
\Delta \psi = - \lambda_d \psi (1-e^{-\psi^2}),
\ee
which is the Euler-Lagrange equation associated with the variational problem
in (\ref{ELuc}). By the results of Section 5c in the same paper, we have $\psi
\in C^2(\R^d)$. Since $\psi$ is $\RSNI$ (recall Lemmas \ref{l2}, \ref{RSNI}),
it follows from (\ref{ELuc}) that $\psi \in L^4(\R^d)$. Suppose that $\psi\in
L^2(\R^d)$. Then, by H\"older's inequality, $\psi \in L^3(\R^d)$. Abbreviate
the right-hand side of (\ref{fax1**}) by $f_{\psi}$. Then $f_{\psi} \in
L^1(\R^d)$,
and so we have
\be{fax1**'}
\psi = f_{\psi} \ast K,
\ee
where $\ast$ denotes convolution and $K$ is the Green function. It follows
that
\be{fax1**''}
\int \psi^2 = \int dy_1 \,f(y_1) \int dy_2 \,f(y_2) \,\Big[ \int dx\,
K(x-y_1) K(x-y_2)\Big].\hskip.5in
\ee
But the last integral is infinite for all $y_1,y_2 \in \R^d$ when $d \leq 4$,
which is a contradiction.
\enddemo

\phantom{goaway}

This proves the existence of a minimiser of (\ref{4.2}) for $2 \leq d \leq 4$.
The remaining claims in Theorem \ref{t4}(i) follow from Lemma \ref{RSNI}.

\demo{{\rm 5.5.} Proof of Theorem {\rm \ref{t4}(ii)}}

\medbreak
  1.  To prove the first claim in Theorem \ref{t4}(ii), we apply
(\ref{scal*}), (\ref{scal**})
with $p=(1-u)^{-1/d}$ and $q=(1-u)^{1/2}$ to (\ref{4.2})  to get
\bel{mon1} 
\quad \chi(u) &=& \inf\left\{\|\nabla \psi\|_2^2 \colon ~\|\psi\|_2=1, 
~\int(e^{-\psi^2}-1+\psi^2) \geq 1-u \right\}   \\[4pt]
&=& (1-u)^{2/d} \inf\Big\{\|\nabla \psi\|_2^2 \colon ~ \|\psi\|_2=1,\nn\\[4pt]
&&\hskip1in
~ \int\Big(\frac{e^{-(1-u)\psi^2}-1+(1-u)\psi^2}{(1-u)^2}\Big) \geq 1\Big\}.
\nne
Since the integrand is nondecreasing in $u$, it follows that the infimum
is nonincreasing in $u$. Hence $u \mapsto (1-u)^{-2/d}\chi(u)$ is
nonincreasing.
To get strict monotonicity we use the existence of a minimiser as
established in
Section 5.4.\medbreak

2.  Let $\psi^*$ be any minimiser of the variational problem in
(\ref{mon1}).
Pick $v \in (u,1)$. Then there exists a $\delta_{u,v}>0$ such that
\be{intineq}
\int\Big(\frac{e^{-(1-v)\psi^{*2}}-1+(1-v)\psi^{*2}}{(1-v)^2}\Big)
\geq 1 + \delta_{u,v}.
\ee
Hence
\bel{ineqchain} 
&&\\
&&\hskip-12pt(1-u)^{-2/d}\chi(u) = \|\nabla\psi^*\|^2_2 \nn\\
&&  \geq \inf\Big\{\|\nabla \psi\|_2^2 \colon ~\|\psi\|_2=1,
~ \int\Big(\frac{e^{-(1-v)\psi^2}-1+(1-v)\psi^2}{(1-v)^2}\Big)
\geq 1 + \delta_{u,v}\Big\}\nn \\
&&  = (1-v)^{-2/d} \inf\Big\{\|\nabla\psi\|^2_2 \colon ~\|\psi\|_2 = 1,\nn\\
&&\hskip1.25in
\int (1-e^{-\psi^2}) \leq 1-(1+\delta_{u,v})(1-v)\Big\},
\nne
where we reverse the scaling that led from (\ref{4.2}) to (\ref{mon1}). But, the
right-hand side is equal to $(1-v)^{-2/d}\chi(v-\delta_{u,v}(1-v))$, which is
strictly larger than $(1-v)^{-2/d}\chi(v)$ by Theorem \ref{t3}(ii). This proves
the strict monotonicity.\medbreak

3.  The second claim in Theorem \ref{t4}(ii) is proved by deriving upper
and lower bounds. Since $e^{-x}\le1-x+\frac{1}{2}x^2$ for $x \geq 0$, we have
\be{4.33}
\int(e^{-\psi^2}-1+\psi^2)\le \frac{1}{2}\int_{\R^d}\psi^{4}.
\ee
Hence (\ref{4.2}) gives
\be{4.34}
\chi(1-\delta)\ge\inf\biggl\{\|\nabla\psi\|_2^2 \colon
~\|\psi\|_2=1,\,\int \psi^4\ge2\delta\biggr\}.
\ee
We use (\ref{scal*}), (\ref{scal**}) with $p=(2\delta)^{-1/d}$ and $q=1$ in
(\ref{4.34}),
to obtain
\be{4.35}
\chi(1-\delta)\ge(2\delta)^{2/d}\inf\{\,\|\nabla\psi\|_2^2 \colon
\|\psi\|_2=1,\,\|\psi\|_4\ge1\,\}
= (2\delta)^{2/d}\mu_d.\quad
\ee
The fact that the infimum equals $\mu_d$ in (\ref{eq1.14}) follows from the
same type
of argument as in the proof of Lemma \ref{l6}, showing that the constraint
$\|\psi\|_4
\geq 1$ may be replaced by $\|\psi\|_4=1$. For the case $d=4$ we shall see
in Step 4
in Section 5.5 that the constraint $\|\psi\|_4=1$ can even be dropped.
Hence we have
$\delta^{-2/d}\chi(1-\delta) \geq 2^{2/d}\mu_d$, which proves the lower bound.\medbreak

4.  Since
$e^{-\psi^2}-1+\psi^2\geq\frac{1}{2}\psi^4-\frac{1}{6}\psi^6$, we have by (\ref{4.2}),
\be{4.43}
\chi(1-\delta)\le\inf\{\|\nabla\psi\|_2^2 \colon ~\|\psi\|_2=1,\,
\|\psi\|^4_4 \geq 2\delta  + \textstyle\frac{1}{3} \|\psi\|^6_6\}.\hskip.5in
\ee
We apply (\ref{scal*}), (\ref{scal**}) with $p=(2\delta)^{-1/d}$ and
$q=(2\delta)^{1/2}$,
to get
\be{fext}
(2\delta)^{-2/d}\chi(1-\delta) \leq \inf\{\|\nabla\psi\|_2^2 \colon
~\|\psi\|_2=1,\,
\|\psi\|^4_4 \geq 1 + \textstyle\frac{2\delta}{3} \|\psi\|^6_6\}.\hskip.25in
\ee

\medbreak
\underline{$d=2$}: Insert the Sobolev inequality $\|\psi\|_6 \leq S_{2,6}
(1 + \|\nabla\psi\|^2_2)^{1/2}$ (Lieb and Loss \cite[p.\ 190]{LL}) into
(\ref{fext}), to get
\be{ld=2'}
(2\delta)^{-1}\chi(1-\delta) \leq \inf\{\|\nabla\psi\|^2_2 \colon
~\|\psi\|_2=1,\,
\|\psi\|^4_4 \geq 1 + \textstyle{\frac{2\delta}{3}}
S^6_{2,6}(1+\|\nabla\psi\|^2_2)^3\}.
\ee
By considering the trial function $\psi_a(x) = (\pi a^2)^{-1/2}
\exp[-|x|^2/2a^2]$
($x \in \R^2$) for $a>0$ sufficiently small, we see that the second
constraint in
(\ref{ld=2'}) is satisfied for $\delta>0$ sufficiently small. Hence
$\|\nabla\psi\|^2_2
\leq M$ for some $M<\infty$ and all $\psi$ satisfying the constraints in
(\ref{ld=2'}).
Thus, for all $\delta>0$ sufficiently small,
\be{ld=2''}
\ba{lll}
(2\delta)^{-1}\chi(1-\delta) &\leq& \inf\{\|\nabla\psi\|^2_2 \colon
~\|\psi\|_2=1,\,
\|\psi\|^4_4 \geq 1 + \textstyle{\frac{2\delta}{3}} S^6_{2,6}(1+M)^3\}\\ \\
&\leq& (1 + \textstyle{\frac{2\delta}{3}} S^6_{2,6}(1+M)^3) \mu_2,
\ea
\ee
where we use (\ref{scal*}), (\ref{scal**}) once more, this time with $p=(1 +
\frac{2\delta}{3}
S^6_{2,6}(1+M)^3)^{1/2}$ and $q=p^{-1}$. We also recall (\ref{eq1.14}).
Let $\delta \da 0$
to get the desired upper bound.

\medbreak
\underline{$d=3$}: Insert the Sobolev inequality $\|\psi\|_6 \leq
S_3^{-1/2}\|\nabla\psi\|_2$ (recall (\ref{4.38}) and (\ref{4.40})) into
(\ref{fext}), to get
\be{ld=3'}
(2\delta)^{-2/3}\chi(1-\delta) \leq \inf\{\|\nabla\psi\|^2_2 \colon
~\|\psi\|_2=1,\,
\|\psi\|^4_4 \geq 1 + \textstyle{\frac{2\delta \|\nabla\psi\|^6_2}{3S_3^3}}\}.\hskip.25in
\ee
This time we use the trial function $\psi_a(x) = (\pi a^2)^{-3/4}
\exp[-|x|^2/2a^2]$
($x \in \R^3$), to see that the second constraint in (\ref{ld=3'}) is
satisfied for
$\delta>0$ sufficiently small, and that $\|\nabla\psi\|^2_2 \leq M$ for
some $M<\infty$ and all
$\psi$ satisfying the constraints in (\ref{ld=3'}). After scaling we get
\be{ld=3''}
(2\delta)^{-2/3}\chi(1-\delta) \leq (1 +
\textstyle{\frac{2\delta M^3}{3S_3^3}})^{2/3}\mu_3,
\ee
and hence the desired upper bound after letting $\delta \da 0$.

\medbreak
\underline{$d=4$}: This case is slightly more subtle.\\

i.  For any $M>0$ we may insert the
constraint $\|\psi\|_{\infty}\leq M$, to get
\be{ld=4a}
(2\delta)^{-1/2}\chi(1-\delta) \leq \inf\{\|\nabla\psi\|^2_2 \colon
~\|\psi\|_2=1,\,
\|\psi\|_{\infty}\leq M,\, \|\psi\|^4_4 \geq 1 +\textstyle\frac{2\delta
M^2}{3}\|\psi\|^4_4\}.
\ee
For any $0<\delta<3/2M^2$ we therefore have
\bel{ld=4b}&&\\
&&\hskip-.35in(2\delta)^{-1/2}\chi(1-\delta)
\nn\\
&& \leq \inf\{\|\nabla\psi\|^2_2 \colon
~\|\psi\|_2=1,\,
\|\psi\|_{\infty}\leq M,\, \|\psi\|^4_4 \geq \left(1 - \frac{2\delta
M^2}{3}\right)^{-1}\} \nn \\
&& = \left(1 - \frac{2\delta M^2}{3}\right)^{-1/2} \inf\{\|\nabla\psi\|^2_2
\colon ~\|\psi\|_2=1,\nn\\
&&\hskip2.075in
\|\psi\|_{\infty} \leq M\left(1 - \frac{2\delta M^2}{3}\right)^{1/2},
\, \|\psi\|^4_4 \geq 1\}.
\nne
We will construct a sequence $(\psi_j)$ in $D^1(\R^4)$ satisfying the
constraints
in (\ref{ld=4b}), with $\|\psi_j\|_4=1$ for all $j$, such that $\limsup_{j
\ra \infty}
\|\nabla\psi_j\|^2_2 \leq S_4$. Since $\mu_4 \geq S_4$, this sequence will be a
minimising sequence of $\mu_4$ (recall (\ref{eq1.14}) and (\ref{4.38})).\medbreak

ii.  Let $\psi_0$ be defined by
\be{ld=4c}
\psi_0 (x) = \left(|x|^2 + \sqrt{\pi^2/6}\right)^{-1} \qquad (x \in \R^4).
\ee
We have $\|\psi_0\|_4=1$. For $\alpha>0$, let $\psi_{\alpha}$ be defined by
\be{ld=4d}
\psi_{\alpha}(x) = c_{\alpha} e^{-\alpha|x|}\psi_0(x) \qquad (x \in \R^4),
\ee
where $c_{\alpha}>0$ is chosen such that $\|\psi_{\alpha}\|_4 = 1$.
Consider now
the scaling
\be{ld=4e}
\psi_{\alpha,\beta}(x) = \beta^{-1} \psi_{\alpha}(x/\beta) \qquad (x \in \R^4)
\ee
with $\beta = \beta(\alpha) > 0$ chosen such that $\|\psi_{\alpha,\beta}\|_2=1$, i.e.,
$\beta = \|\psi_a\|_2^{-1}$. We have $\|\psi_{\alpha,\beta}\|_4 = 1$,
while
\be{ld=4f}
\|\psi_{\alpha,\beta}\|_{\infty} = c_{\alpha} \beta^{-1} \|\psi_0\|_{\infty}
= c_{\alpha} \|\psi_{\alpha}\|_2 \sqrt{6/\pi^2}
\ee
with
\be{ld=4g}
\|\psi_{\alpha}\|^2_2 = c_{\alpha}^2 \int e^{-2\alpha |x|} \psi_0^2(x) dx
\leq c_{\alpha}^2 \int_0^{\infty} e^{-2\alpha r} r^{-7/2} [2\pi^2 r^3] dr
\leq C c_{\alpha}^2 \alpha^{-1/2}.
\ee
Hence $\|\psi_{\alpha,\beta}\|_{\infty} \leq C' c_{\alpha}^2 \alpha^{-1/4}$.
Now we pick $\alpha = \delta$ and $M = \delta^{-1/3}$, and we\break note that $C'
c_{\delta}^2
\delta^{-1/4} \leq \delta^{-1/3}(1-\frac{2\delta^{1/3}}{3})^{1/2}$ for $\delta$
sufficiently small because $\lim_{\delta \da 0} c_{\delta} = 1$. Then
$\psi_{\alpha,\beta}$
satisfies the constraints in (\ref{ld=4b}), and so for $\delta$
sufficiently small,
\be{ld=4g*}
(2\delta)^{-1/2}\chi(1-\delta) \leq (1 -
\textstyle\frac{2\delta^{1/3}}{3})^{-1/2}
\|\nabla\psi_{\alpha,\beta}\|^2_2 \quad \hbox{with } \alpha = \delta,
\beta = \beta(\delta).\hskip12pt
\ee
\bigbreak iii.  It thus remains to show that $\psi_{\alpha,\beta}$ is a
minimising sequence
for the variational problem defining $\mu_4$. To that end we first note that
$\|\nabla\psi_0\|^2_2 = S_4$. Therefore, using (\ref{ld=4d}) and (\ref{ld=4g}),
we have
\bel{ld=4h}
 \|\nabla\psi_{\alpha,\beta}\|^2_2 &
  = &\|\nabla\psi_{\alpha}\|^2_2\\[4pt]
&= &\alpha^2 \|\psi_{\alpha}\|^2_2
+ c_{\alpha}^2 \int e^{-2\alpha |x|} |\nabla\psi_0(x)|^2 dx\nn\\[4pt]
&&
+ 2\alpha c_{\alpha} \int \psi_{\alpha}(x) e^{-\alpha |x|}
|\nabla\psi_0(x)| dx\nn \\[4pt]
& \leq &\alpha^2 \|\psi_{\alpha}\|^2_2
+  c_{\alpha}^2 \|\nabla\psi_0\|^2_2
+ 2\alpha c_{\alpha} \|\psi_{\alpha}\|_2 ~\|\nabla\psi_0\|_2\nn\\[4pt]
&=& \alpha^2 \|\psi_{\alpha}\|^2_2
+ c_{\alpha}^2 S_4
+ 2\alpha c_{\alpha} \|\psi_{\alpha}\|_2 S_4^{1/2}\nn \\[4pt]
&  \leq& C c_{\alpha}^2 \alpha^{3/2}
+ c_{\alpha}^2 S_4
+ 2 C^{1/2} c_{\alpha}^2 \alpha^{3/4} S_4^{1/2}.
\nne
Since $\lim_{\alpha \da 0} c_{\alpha} = 1$, we can combine this with
(\ref{ld=4g*}) to arrive at
\be{ld=4i}
\limsup_{\alpha \da 0} \|\nabla\psi_{\alpha,\beta}\|^2_2 \leq S_4.
\ee
Finally we note that $\mu_4 \geq S_4$, giving
\be{ld=4j}
\limsup_{\delta \da 0} (2\delta)^{-1/2}\chi(1-\delta) \leq \mu_4.
\ee

\phantom{go away}

  5.  We complete this section by proving that $\mu_d>0$ for $2 \leq d
\leq 4$.
\specialnumber{15} \proclaim{Lemma} \label{r4.4}
\be{4.56}
\mu_2 \ge 1/4S_{2,4}^4 = 27\pi/16,\quad
\mu_3 \ge S_3 = 3(\pi/2)^{4/3},\quad
\mu_4 = S_4 = 4\pi\sqrt{6}/3.
\ee
\el

\vglue9pt
\demo{Proof}
(1) Rewrite (\ref{4.36}) as follows:
\be{4.57}
\|\nabla f\|_2^2 \geq S_{2,4} \|f\|_4^2-\|f\|_2^2.
\ee
Using the scaling in (\ref{scal*}), (\ref{scal**}) with $p>0$ arbitrary and
$q=1$, we have
\be{4.58}
\|\nabla f\|_2^2\geq S_{2,4}p\|f\|_4^2-p^2\|f\|_2^2.
\ee
Putting $\|f\|_2=\|f\|_4=1$ and optimizing over $p$, we get the bound for
$\mu_2$.
\bigbreak
(2) The bound for $\mu_3$ follows from (\ref{eq1.14}) via
(\ref{4.38}), (\ref{4.40}).
\bigbreak
(3) The identity for $\mu_4$ was already proved in Step 4.
\enddemo
 
\demo{{\rm 5.6.} Proof of Theorem {\rm \ref{t5}(i)}}
This section has three parts:\begin{itemize}
\item[(I)] Proof of $0<\nu_d<\infty$. 
\item[(II)] Proof of $\Sigma \neq \emptyset$. 
\item[(III)] Proof of (\ref{eq1.16}). \end{itemize}

(I) Proof of $0<\nu_d<\infty$.
The upper bound is easily deduced from (\ref{eq1.15}) by substituting a
test function. The lower bound comes from the following.

\specialnumber{1}\numbereddemo{{R}emark}\label{r4.6}
{\it For} $d\ge5$
\be{4.67}
\nu_d\ge2^{(d-2)/d}S_d,
\ee
{\it where $S_d$ is the sharp constant in the Sobolev inequality given by}
(\ref{4.40}).
\enddemo
\bigbreak
{\it Proof}.
Since $e^{-x}-1+x \leq \frac12x^\alpha$ for $x\geq0$ and $1\le\alpha\le 2$,
it follows that (recall (\ref{eq1.15}))
\be{4.68}
\nu_d\ge\inf\{\, \|\nabla\psi\|_2^2 \colon ~
\textstyle\int \textstyle\frac12 \psi^{2\alpha} \geq 1\,\}
=2^{1/\alpha}\inf\{\,\|\nabla\psi\|_2^2 \colon ~\|\psi\|_{2\alpha}\ge 1\,\}.\hskip.25in
\ee
Pick $\alpha=d/(d-2)$. Then (\ref{4.67}) follows from (\ref{4.68}) via
(\ref{4.38}), (\ref{4.40}).
\hfill\qed \bigbreak
(II) Proof of  $\Sigma \neq \emptyset$.
\medbreak
  1.  We begin with a tail estimate.
\specialnumber{16} \proclaim{Lemma} \label{ldot}
Let $d \geq 3$ and let $\psi \in D^1(\R^d)$ be $\RSNI${\rm .} Then
\be{fax1'}
0 \le \psi(r)\le\biggl(\frac{C}{S_d}\biggr)^{1/2}
\omega_d^{-(d-2)/2d}r^{-(d-2)/2}
\ee
with $\omega_d=|B_1(0)|${\rm ,} $r=|x|$ and $C=\|\nabla \psi\|_2^2${\rm .}
\el

\demo{Proof}
By the Sobolev inequality (\ref{4.38}), (\ref{4.40}), we have
\bel{fax3'}
\|\nabla \psi\|_2^2  &\ge&  S_d \|\psi\|_{2d/(d-2)}^2
\ge S_d \|\psi 1_{B_r(0)}\|_{2d/(d-2)}^2 \\[4pt]
&\ge& S_d \psi(r)^2 |B_r(0)|^{(d-2)/d}
= S_d \psi(r)^2 \omega_d^{(d-2)/d} r^{d-2}.
\nne
\vglue-20pt\enddemo

\vglue12pt

  2. Let $(\psi_j)$ be a minimising sequence for the variational
problem in (\ref{eq1.15}).
We may assume that $\psi_j$ is $\RSNI$ (recall Section 5.1), and that
$\|\nabla \psi_j\|_2^2
\da \nu_d$. We can extract a subsequence, again denoted by $(\psi_j)$, such
that $\psi_j \to
\psi^*$ weakly in $D^1(\R^d)$ and $\psi_j \to \psi^*$ almost everywhere as
$j \to \infty$.
It follows that $\psi^*$ is $\RSNI$ too. Moreover, $\nu_d \geq \|\nabla
\psi^* \|_2^2$.
It therefore suffices to show that $\psi^*$ satisfies the constraint in
(\ref{eq1.15}), since
this implies that $\nu_d \leq \|\nabla \psi^* \|_2^2$, and hence that
$\psi^*$ is a minimiser.\medbreak

  3. Let $\varepsilon >0$ be arbitrary. Since $d \geq 5$, we have $0 \leq
e^{-x}-1+x \leq
x^{d/(d-2)}$ for $x \geq 0$. Hence (\ref{4.38}) and (\ref{4.40}) give
\be{fax6'}
0 \leq \int (e^{-\psi^{*2}}-1+\psi^{*2}) \leq \int \psi^{*2d/(d-2)}
\leq \Big(\frac{\nu_d}{S_d}\Big)^{(d-2)/d}\hskip.5in
\ee
and so there exists an $R_1(\varepsilon)$ such that
\be{fax7'}
0 \leq \int_{B^c_{R_1(\varepsilon)}}(e^{-\psi^{*2}}-1+\psi^{*2}) \leq \varepsilon.
\ee
Let $C=\sup_j \|\nabla\psi_j\|_2^2<\infty$ and define $R_2(\varepsilon)$ by
\be{fax8'}
\Big(\frac{C}{S_d}\Big)^2 \omega_d^{-(d-4)/d} \frac{d}{2(d-4)}
R_2(\varepsilon)^{-(d-4)}
= \varepsilon.
\ee
Then with the help of Lemma \ref{ldot} we obtain
\be{fax9'}
\int\limits_{B^c_{R_2(\varepsilon)}}(e^{-\psi_j^2}-1+\psi_j^2)
\leq \int\limits_{B^c_{R_2(\varepsilon)}} \frac{1}{2}\psi_j^4
\leq \Big(\frac{C}{S_d}\Big)^2 \omega_d^{-2(d-2)/d}
\frac{1}{2} \int\limits_{B^c_{R_2(\varepsilon)}} |x|^{4-2d}dx
= \varepsilon.
\ee

  4.  Put $R(\varepsilon) = \max\{R_1(\varepsilon),R_2(\varepsilon)\}$. Since
$\int (e^{-\psi_j^2}
-1+\psi_j^2)=1$ for all $j$, we get from (\ref{fax9'}) that
\be{fax11'}
1-\varepsilon \leq \int_{B_{R(\varepsilon)}}
(e^{-\psi_j^2}-1+\psi_j^2) \leq 1.
\ee
Moreover, by Lemma \ref{ldot},
\be{fax13'}
0 \leq e^{-\psi_j^2}-1+\psi_j^2 \leq \psi_j^2 \leq \frac{C}{S_d}
\omega_d^{-(d-2)/d} r^{-(d-2)},
\ee
where the right-hand side is integrable on $B_{R(\varepsilon)}$. Since $\psi_j
\to \psi^*$
almost everywhere as $j \to \infty$, we therefore have by the dominated
convergence
theorem and (\ref{fax11'}) that
\be{fax14'}
1-\varepsilon \leq \int_{B_{R(\varepsilon)}}(e^{-\psi^{*2}}-1+\psi^{*2})
\leq 1.
\ee
Combining (\ref{fax7'}) and (\ref{fax14'}), we obtain
\be{fax15'}
1-\varepsilon \leq \int (e^{-\psi^{*2}}-1+\psi^{*2}) \leq 1+\varepsilon.
\ee
Since $\varepsilon$ was arbitrary, we conclude that $\psi^*$ satisfies the
constraint in (\ref{eq1.15}). By the same argument as in the proof of Lemma
\ref{RSNI},
we get that $\psi^*$ is strictly decreasing in the radial component.

\medbreak
(III) Proof of (\ref{eq1.16}).
\smallbreak
 1.  We begin with the lower bound. Let $\psi$ be any minimiser of the
variational
problem (\ref{eq1.15}) centered at 0. By the results of Section 5b in
Berestycki and Lions
\cite{BeLi}, there exists a Lagrange multiplier $\lambda_d>0$ such that
\be{fax1}
(r^{d-1}\psi')' = - \lambda_d r^{d-1} \psi (1-e^{-\psi^2}),
\ee
which is the Euler-Lagrange equation in the radial form. By the results of
Section 5c
in the same paper, we have $\psi \in C^2(\R^d)$. Because $\psi$ is radially
symmetric and
centered at 0, it follows that $\psi(0)<\infty$ and $\psi'(0)=0$. Hence
$\psi \in
L^{\infty}_{\rm loc}(\R^d)$, and we already know that $\nabla \psi \in
L^2(\R^d)$ and
$\lambda_d(e^{-\psi^2}_{\phantom{|}}-1+\psi^2) \in L^1(\R^d)$. It therefore follows from
Pohozaev's
identity (see Proposition 1 in the same paper) that
\be{fax3}
\|\nabla \psi \|_2^2 = \frac{d}{d-2} \lambda_d \int (e^{-\psi^2}-1+\psi^2)
= \frac{d}{d-2} \lambda_d.
\ee
Hence the Lagrange multiplier can be identified as
\be{fax4}
\lambda_d = \frac{d-2}{d}\nu_d.
\ee
Multiply both sides of (\ref{fax1}) by $\psi$, integrate over $r \in [0,\infty)$
and use integration by parts, to get
\be{fax4ext}
\|\nabla\psi\|_2^2 = \lambda_d \int \psi^2(1-e^{-\psi^2}).
\ee
Combining this with (\ref{fax3}), we obtain
\be{fax4ext'}
\int \psi^2(1-e^{-\psi^2}) = \frac{d}{d-2},
\ee
which obviously implies that $\|\psi\|_2^2>d/(d-2)$.\medbreak

2.  The upper bound is more laborious. We first consider the
case
$d\geq 6$. Then
\be{MB6}
\psi(1-e^{-\psi^2}) \leq \psi^{2d/(d-2)},
\ee
and hence, with the help of (\ref{fax1}) and the Sobolev inequality in
(\ref{4.38}),
we may estimate
\bel{MB7}
\qquad -r^{d-1}\psi'(r) &\leq& \lambda_d \int_0^r s^{d-1} \psi(s)^{2d/(d-2)}ds
\leq \frac{\lambda_d}{d\omega_d} \int_{\R^d} \psi^{2d/(d-2)}  \\[4pt]
&\leq&
\frac{\lambda_d}{d\omega_d}\Big(\frac{\|\nabla\psi\|^2_2}{S_d}\Big)^{d/(d-2)
}
= \frac{\lambda_d}{d\omega_d} \Big(\frac{\nu_d}{S_d}\Big)^{d/(d-2)}.
\nne
Multiplying both sides of (\ref{MB7}) by $r^{1-d}$ and integrating over
$[r,\infty)$,
we get
\be{MB8}
\psi(r) \leq A r^{-(d-2)} \mbox{ with } A = \frac{\lambda_d}{d(d-2)\omega_d}
\Big(\frac{\nu_d}{S_d}\Big)^{d/(d-2)}.\hskip.5in
\ee
Next, by (\ref{MB8}) and the Sobolev inequality,
\bel{MB9}
&&\\
\int \psi^2 &=& \int_{\{\psi\geq 1\}} \psi^2 + \int_{\{\psi<1\}} \psi^2
\leq \int_{\{\psi\geq 1\}} \psi^{2d/(d-2)} + \int_{\{\psi<1\}}
\min\{\psi^2,1\}  \nn\\[4pt]
&\leq& \left(\frac{\nu_d}{S_d}\right)^{d/(d-2)} + d\omega_d
\int_0^{\infty} \min\{A^2 r^{4-2d},1\} r^{d-1} dr,
\nne
which is finite because $d \geq 6$.
\medbreak
  3.  We note that (\ref{MB6}) fails for $d=5$. But
$\psi(1-e^{-\psi^2})\leq \psi^3$,
and so because $\psi$ is $\RSNI$ we have, by the Sobolev inequality and
(\ref{fax1}),
\be{MB11}
-r^4\psi'(r) \leq \lambda_5 \int_0^r s^4 \psi(s)^3 ds
\leq \psi(r)^{-1/3} \frac{\lambda_5}{5\omega_5} \int \psi^{10/3}
\leq  \psi(r)^{-1/3} \frac{\lambda_5}{5\omega_5} \left(\frac{\nu_5}{S_5}\right)^{5/3}.
\ee
Hence $(\psi(r)^{4/3})' \geq -C_1 r^{-4}$. Integrating this inequality over
$[r,\infty)$,
we find $\psi(r) \leq C_2 r^{-9/4}$. Returning to (\ref{fax1}) once more,
we have by
H\"older's\break inequality,
\bel{MB13}
&&\\
-r^4 \psi'(r)& \leq &\lambda_5 \int_0^r s^4 \psi(s)^3 ds
\nn\\
&\leq& \lambda_5 \Big(\int_0^r s^4 \psi(s)^{10/3}ds\Big)^{4/5}
\Big(\int_0^r s^4 \psi(s)^{5/3} ds\Big)^{1/5}
\leq C_3 r^{1/4},\nne
where we use the Sobolev inequality for the first integral and our bound on
$\psi(r)$
for the second integral. Multiplying both sides of (\ref{MB13}) by $r^{-4}$
and integrating
over $[r,\infty)$, we arrive at $\psi(r) \leq C_4 r^{-11/4}$. But this
bound is integrable
on the set where $C_4 r^{-11/4}<1$, and so we obtain that $\psi \in
L^2(\R^5)$ via an
estimate similar to (\ref{MB9}).

\demo{{\rm 5.7.} Proof of Theorem {\rm \ref{t5}(ii)}}
By Lemma \ref{l6} and (\ref{scal*}), (\ref{scal**}) with $p=(1-u)^{1/d}$ and
$q=1$,
\be{eq4.151}
\chi(u) = \widetilde \chi(u)
= (1-u)^{(d-2)/d}\inf\{\|\nabla\psi\|_2^2 \colon ~\|\psi\|_2^2 \leq
(1-u)^{-1},\textstyle\int(e^{-\psi^2}-1+\psi^2)=1\}.
\ee

 1.  Let $u \in (0,u^*_d]$. There exists a $\RSNI$ (modulo shifts)
minimising
sequence $(\psi_j)$ of the variational problem for $\chi(u)$ such that
$\|\nabla\psi_j\|^2_2 \ra \chi(u)$ and $\psi_j \ra \psi^*$ weakly in
$H^1(\R^d)$ as
$j \to \infty$. By extracting a subsequence, again denoted by $(\psi_j)$,
we also have
$\psi_j \ra \psi^*$ almost everywhere as $j \to \infty$. Suppose that
$\|\psi^*\|^2_2
< (1-u)^{-1}$. Then we can do smooth perturbations of $\psi^*$ inside the class
$\{\psi \in H^1(\R^d) \colon ~\int (e^{-\psi^2}-1+\psi^2)=1\}$ to conclude
that $\psi^*$ must be an element of $\Sigma^*$, the set of local minimisers of
the variational problem in (\ref{eq1.15}). But by the definition of $u^*_d$
we have
$\|\psi^*\|^2_2 \geq (1-u^*_d)^{-1} \geq (1-u)^{-1}$, which is a contradiction.
Therefore $\|\psi^*\|^2_2 = (1-u)^{-1}$, and so $\psi^*$ is a minimiser of
$\chi(u)$.
Clearly, $\psi^*$ is $\RSNI$ (modulo shifts). By Lemma \ref{RSNI} it must
be strictly
positive and strictly decreasing in the radial component.\medbreak

  2.  Let $u \in (u^+_d,1)$. Suppose that the variational problem for
$\chi(u)$ has a
minimiser $\psi^*$. Then $\phi^*(\cdot) = \psi^*(\cdot\,(1-u)^{1/d})$ is a
minimiser of the
variational problem in (\ref{eq4.151}). Since $u>u^+_d\geq u^-_d$, it
follows from (\ref{eq1.17})
that $\|\nabla\phi^*\|_2^2 =~\nu_d$. Hence $\phi^*$ is also a minimiser of
(\ref{eq1.15}).
But all such minimisers have a squared $L^2$--norm that is bounded above by
$(1-u^+_d)^{-1}$.
This is a contradiction because $\|\phi^*\|_2^2 =
(1-u)^{-1}\|\psi^*\|^2_2=(1-u)^{-1}>
(1-u^+_d)^{-1}$. Hence $\chi(u)$ has no minimiser. The last claim in
Theorem \ref{t5}(ii)
is obvious.\enddemo

\demo{{\rm 5.8.} Proof of Theorem {\rm \ref{t5}(iii)}}
\smallbreak
 1.  By dropping the constraint $\|\psi\|_2^2 \leq (1-u)^{-1}$ from
(\ref{eq4.151})
and recalling (\ref{eq1.15}), we see that $\chi(u)\geq
(1-u)^{(d-2)/d}\nu_d^{\phantom{|}}$, which proves
the lower bound in (\ref{eq1.17}). The upper bound is proved via an
argument similar
to Step ii in the proof of Lemma \ref{l6}. Let $u \in (u^-_d,1)$. Then, by
Theorem \ref{t5}(i), there
exists a $\psi^* \in \Sigma$ satisfying $\|\psi^*\|^2_2<(1-u)^{-1}$. For $n
\in \N$,
define $\psi^*_{n,u}$ by
\be{eq4.147}
\psi^{*2}_{n,u} = \psi^{*2} + \Big(\frac{1}{1-u} - \|\psi^*\|_2^2\Big) ~p_n,
\ee
with $p_n$ as in (\ref{pdef1}). Then $\|\psi^*_{n,u}\|_2^2=(1-u)^{-1}$ for
all $n$. Moreover,
since $x\mapsto e^{-x}-1+x$ is increasing on $[0,\infty)$, we have
\be{eq4.149}
\int(e^{-\psi^{*2}_{n,u}}-1+\psi^{*2}_{n,u})\ge
\int(e^{-\psi^{*2}}-1+\psi^{*2})=1.
\ee
Therefore, by (\ref{eq4.151}) and the convexity inequality for gradients,
\bel{eq4.150}
\qquad(1-u)^{-(d-2)/d}\chi(u)
& \le& \|\nabla\psi^*_{n,u}\|_2^2\\
&\le& \|\nabla\psi^*\|_2^2 + \Big(\frac{1}{1-u} - \|\psi^*\|_2^2\Big)
|\sqrt{\nabla p_n}\|_2^2\nn\\
&\le& \nu_d + (1-u)^{-1}\frac{2d}{n^2}.
\nne
Let $n\to\infty$ to get the upper bound in (\ref{eq1.17}). The case
$u=u^-_d$ follows by
continuity of $\chi$.

  2. It is immediate from (\ref{eq4.151}) that $u \mapsto (1-u)^{-
(d-2)/d}\chi(u)$ is
nonincreasing on $(0,1)$. Strict monotonicity on $(0,u^*_d)$ follows from
the existence
of a minimiser via the same type of argument as in Step 2 in Section 5.5.
Finally, the
fact that $u \mapsto (1-u)^{-(d-2)/d}\chi(u) > \nu_d$ for all $u \in
(u^*_d,u^-_d)$ needs
to be proved only when $u^*_d<u^-_d$. In that case $\Sigma^*\setminus\Sigma
\neq \emptyset$.
But clearly $\|\nabla\psi\|^2_2>\nu_d$ for all
$\psi\in\Sigma^*\setminus\Sigma$, and so
weak convergence to any such $\psi$ will not reach the minimal value $\nu_d$.

\medbreak
We conclude by settling an old debt: to prove the strict monotonicity in
Theorem
\ref{t3}(iii). By Theorems \ref{t4}(i) and \ref{t5}(ii), (\ref{chiscal3}) has a
minimiser for $2\le d \le 4$, $u\in (0,1)$, and for $d \ge 5$, $u \in (0,u^*_d]$,
in which case the claim follows via the same type of argument as in Step 2 in
Section 5.5. On the other hand, for $d \ge 5$, $u \in (u^*_d,1)$, we can
appeal to Theorem \ref{t5}(iii), which easily gives the claim because
$u^*_d \geq 2/d$.

\bye